\documentclass[11pt,fleqn]{article}
\usepackage{latexsym,amsmath,amssymb}
\usepackage{longtable}
\usepackage{hyperref}
\usepackage{array}
\usepackage{graphics}
\usepackage{multirow}
\usepackage{a4wide}
\usepackage{epsfig}
\usepackage{graphicx}
\usepackage{caption}
\usepackage{subcaption}
\usepackage[title]{appendix}
\usepackage{url}
\usepackage{float}
\usepackage{bm}
\usepackage{xcolor}
\begin{document}
\numberwithin{equation}{section}
\renewcommand{\thefootnote}{*}
\thispagestyle{empty}

\begin{center}
{\Large \textbf { \vspace{1.5mm} Lower $k-$record values from unit-Gompertz distribution and associated inference}}\\
\vspace{.5cm}
{\bf}~{\bf Ehsan Ormoz$^1$, Zuber Akhter$^2$\footnote{Corresponding author e-mail address: akhterzuber022@gmail.com}, Mahfooz Alam$^3$ and \small\bf S.M.T.K. MirMostafae$^4$ }\\

\vspace{.3cm}
$^1$Department of Statistics\\ Mashhad Branch, Islamic Azad University, Mashhad, Iran

\vspace{.3cm}
$^2$Department of Statistics\\
University of Delhi, Delhi-110 007, India

\vspace{.3cm}
$^3$Department of Mathematics $\&$ Statistics\\
Vishwakarma University, Pune 411 048, Maharashtra, India

\vspace{.3cm}
$^4$Department of Statistics\\
University of Mazandaran, P.O. Box 47416-1467, Babolsar, Iran

\end{center}

\vspace{0.25cm}
{\bf Abstract.} Mazucheli $et~al.$ (2019) introduced the unit-Gompertz (UG) distribution and studied some of its properties. More specifically, they considered the random variable $X = {\rm e}^{-Y}$, where $Y$ has the Gompertz distribution. In this paper, we consider the lower $k$-record values from this distribution. We obtain exact explicit expressions as well as several recurrence relations for the single and product moments of lower $k$-record values and then we use these results to compute the means, variances and the covariances of the lower $k$-record values. We make use of these calculated moments to find the best linear unbiased estimators (BLUEs) of the location and scale parameters of the UG distribution. Applying the relation between the BLUE and the best linear invariant estimator (BLIE), we obtain the BLIEs of the location and scale parameters, as well. In addition, based on the observed $k$-records, we investigate how to obtain the best linear unbiased predictor (BLUP) and best linear invariant predictor (BLIP) for a future $k$-record value. Confidence intervals for the unknown parameters and prediction intervals for future $k$-records are also discussed.
A simulation study is performed to assess the point and interval estimators and predictors proposed in the paper. The results show that the BLIE and BLIP outperform the BLUE and BLIP, in the sense of mean squared error criterion, respectively.   Finally, a real data set pertaining to COVID-19 2-records is analyzed.

\vspace{.25cm}

\section{Introduction}\label{sec1}
Mazucheli $\it{et \, al.}$ (2019) proposed the unit-Gompertz (UG) distribution by the transformation of $X={\rm e}^{-Y}$, a new distribution hold on (0,1), where $Y$ follows the Gompertz distribution. The probability density function (pdf) and cumulative distribution function (cdf) of the UG distribution are
\begin{equation*}\label{PDF1}
f(x)=\theta\,\alpha\,x^{-{(\theta+1)}}\,\exp\Big[-\alpha\big({x^{-\theta}-1}\big)\Big]
\end{equation*}
and
\begin{equation*}
F(x)=\exp\Big[-\alpha\big({x^{-\theta}-1}\big)\Big],
\end{equation*}
respectively, where $x\in(0,1),\, \alpha>0 $ and $\theta>0$.

\vspace{0.25cm}
Now, let us introduce the three parameter UG distribution  (the location parameter is added), say UG$(\mu,\sigma,\theta)$, by
\begin{equation}\label{pdf}
f(x)=\frac{\theta}{\sigma}\,\left(\frac{x-\mu}{\sigma}\right)^{-{(\theta+1)}}\,\exp\Big[-\Big({\left\{\frac{x-\mu}{\sigma}\right\}^{-\theta}-1}\Big)\Big],
\end{equation}
where  $x\in(\mu,\mu+\sigma),\, \mu\in \mathbb{R},\, \sigma>0$ and $\theta>0$. The parameters $\mu,~\sigma$ and $\theta$ are the location, scale and shape parameters, respectively.

\vspace{0.25cm}
Note that the standardized UG distribution has the pdf

\begin{equation}\label{pdfs}
f(z)={\theta}{z}^{-{(\theta+1)}}\,\exp\Big[-\left({{z}^{-\theta}-1}\right)\Big],  \quad z\in(0,1),\, \theta>0,
\end{equation}

and cdf

\begin{equation}\label{cdfs}
F(z)=\exp\Big[-\left({{z}^{-\theta}-1}\right)\Big],  \quad z\in(0,1).
\end{equation}

Henceforth, we denote by $Z$ a random variable with density (\ref{pdfs}), i.e. $Z\sim$UG$(0,1,\theta)$. Notice that
\begin{equation}\label{2.3}
f(z)=\theta\,z^{-(\theta+1)}F(z).
\end{equation}
Mazucheli $\it{et \, al.}$ (2019) showed in many ways that a unit-Gompertz distribution is a better model than some well-known unit distributions, such as Beta and Kumaraswamy distributions. A unit-Gompertz distribution exhibits constant, increasing, upside-down bathtub shaped hazard rate modeling. It can be a suitable model for fitting the skewed data which may not be properly fitted by other common distributions and can also be used in a variety of problems in different areas such as economic, agriculture, and marketing industry areas as only the limited range such as percentages, proportions or fractions. Papke and Wooldridge (1996), Genç (2013) and  Cook $\it{et \, al.}$ (2008) discussed several examples based on proportions/ratios in the unit interval used in empirical finance and reliability theory when the models defined on the unit interval. For more properties of the UG distribution, one may refer to Anis and De (2020), Anis (2020) and Arshad $et~al.$ (2023).

\vspace{0.25cm} 
The concept of record values was originally introduced by Chandler (1952). An observation is referred to as a lower record if its value is less than all the previous observations, while it is considered an upper record if its value is greater than all the preceding observations. Record values find application in a broad spectrum of practical scenarios, including industrial stress testing, meteorological analysis, hydrology, seismology, oil and mining surveys, as well as sports and athletic events. Dziubdziela and  Recently, Kumar $\it{et \, al.}$ (2020) and Ahmed and Aftab (2023) worked on the inferential problems for the UG distribution based on record data.

\vspace{0.25cm}
Kopcoi$\acute{n}$ski (1976) extended Chandler's (1952) concept of record values by introducing random variables with a more generalized nature, which they termed as "$k$-th record values." Subsequently, Minimol and Thomas (2013) also referred to the record values defined by Dziubdziela and Kopcoi$\acute{n}$ski (1976) as "generalized record values." When setting $k=1$, this corresponds to ordinary record statistics.
Kamps (1995) and Danielak and Raqab (2004) emphasized that in certain situations, record values themselves are viewed as outliers, and hence the second or third largest values are of special interest. This issue is avoided once we consider the model of $k$-th record values.

\vspace{0.25cm}
Let $\{X_n, \,n\geq 1\}$ be a sequence of independent and identical distributed  $(iid)$ continuous random variables with cdf $F(x)$ and pdf $f(x)$. The $j$-th order statistic of a sample $X_1, X_2, \dots, X_n$ is denoted by $X_{j:n}$. For a fixed positive integer $k$, we define the sequence $\{L_{k}(n), n\geq1\}$ of $k$-th lower record times of $\{X_n,\, n\geq1\}$ as follows
\begin{equation*}
L_{k}(1)=1,
\end{equation*}
\begin{equation*}
L_{k}(n+1)=\min\{j>L_{k}(n):X_{k:L_{k}(n)+k-1}>X_{k:j+k-1}\}.
\end{equation*}
The sequence $\{Z_{n(k)}, \,n\geq 1\}$ with $Z_{n(k)}=X_{k:L_{k}(n)+k-1}$, $n=1,2,\dots,$ is called the sequence of $k$-th lower record values of $\{X_n, \,n\geq 1\}$. For convenience, we shall also take $Z_{0(k)}=0$. Note that for $k=1$, we have the record values of $\{X_n, \,n\geq 1\}$. For $n\geq1$, the joint density function of the first $n$ $k$-th lower record values is obtained as
\begin{equation}
f_{Z_{1(k)},Z_{2(k)},\dots, Z_{n(k)}}(z_1,z_2,\dots, z_n)=k^n [F(z_n)]^{k}\prod_{i=1}^{n}\frac{f(z_i)}{F(z_i)}.
\end{equation}
The pdf of $Z_{n(k)}$ and the joint pdf of $Z_{m(k)}$ and $Z_{n(k)}$ are respectively as follows (see for example Pawlas and Szynal, 1998)
\begin{equation}
f_{Z_{n(k)}}(z)=\frac{k^n}{(n-1)!} [-\ln {F}(z)]^{n-1}[{F}(z)]^{k-1}f(z),\quad n\geq1,~k\geq1
\end{equation}
and
\begin{equation*}
f_{Z_{m(k)},Z_{n(k)}}(x,y)=\frac{k^n}{(m-1)! \, (n-m-1)!} [-\ln {F}(x)]^{m-1}\frac{f(x)}{{F}(x)}
\end{equation*}
\begin{equation*}
\hspace{1.1in}\times [\ln{F}(x)-\ln{F}(y)]^{n-m-1}[{F}(y)]^{k-1}f(y),
\end{equation*}
\vspace{-0.8cm}
\begin{equation}
\hspace{1.27in} y<x, \, 1\leq m<n, \, n> 1, \, k\geq1.
\end{equation}

Setting $k=1$, the aforementioned configuration of $k$-th lower record values reduces to lower record values.

\vspace{0.25cm}
For some developments in the field of $k$-th lower record values, particularly those originating from distributions such as exponential, inverse Weibull, generalized extreme value, power, uniform, Frech\'et, Gumbel, inverse Pareto, inverse generalized Pareto, inverse Burr and linear exponential distributions, one may refer to the works of Pawlas and Szynal (1998), Pawlas and Szynal (2000a), Pawlas and Szynal (2000b), Bieniek and Szynal (2002), Bieniek and Szynal (2013). MirMostafaee $\it{et \, al.}$ (2016) discussed the estimation and prediction problems for the Topp-Leone distribution based on lower $k$-record values. Recently, Makouei $\it{et \, al.}$ (2021) derived explicit expressions for single and product moments of order statistics and $k$-record values from the complementary beta distribution. These results have subsequently been employed to estimate the parameters of the complementary beta distribution. To validate the practical utility of these theoretical findings, the authors also provided two real data examples.

\vspace{.25cm}
This paper is organized as follows: The single and product moments of lower $k$-records as well as recurrence relations are discussed in Section 2. 
The best linear unbiased estimation and the best linear invariant estimation of the location and scale parameters are developed in Sections 3 and 4, respectively. Besides, making use of the BLUEs and BLIEs, we propose some pivot quantities to construct confidence intervals for the location and scale parameters. In Section 5, we deal with the best linear unbiased prediction and the best linear invariant prediction of a future lower $k$-record value. We also make two pivot quantities to find prediction intervals for a future $k$-record. A simulation study is conducted in Section 6 and a real data set is analyzed in Section 7 to illustrate the theoretical results of the paper and assess the proposed point and interval estimators and predictors. Several conclusions are provided in Section 8, which ends the paper.

\section{Single and Product Moments of Lower $k$-Record Values}\label{sec1}
Let $R_{1(k)},R_{2(k)},\ldots,R_{n(k)}$ be the first $n$ lower \emph{k}-record values arising from the UG distribution with pdf (\ref{pdf}). Further, let $Z_{i(k)}=\dfrac{R_{i(k)}-\mu}{\sigma}$, $i= 1,2,\ldots,n$, be the corresponding  lower \emph{k}-record values, which can be interpreted as the lower $k$-record values from the standardized UG  distribution with pdf (\ref{pdfs}).

\vspace{0.25cm}
Next, the $r$-th single moment of $Z_{n(k)}$ takes the form
\begin{equation}\label{mean}
\mu_{n;k}^{(r)}=E[\big(Z_{n(k)}\big)^r]=\int_{0}^{1}x^{r}f_{Z_{n(k)}}(x){\rm d}x,\quad n\geq 1,~r\in\mathbb{N},
\end{equation}

and the $(r,s)$-th product moment of $Z_{n(k)}$ and $Z_{m(k)}$ reduces to
\begin{equation}\label{pmean}
\mu_{m,n;k}^{(r,s)}=E[\big(Z_{m(k)}\big)^r\big(Z_{n(k)}\big)^s]=\int_{0}^{1}\int_{y}^1 x^{r}y^{s}f_{Z_{m(k)},Z_{n(k)}}(x,y){\rm d}x{\rm d}y,\quad n>m\geq 1,~r,s\in\mathbb{N}.
\end{equation}

\subsection{Single moments}\label{sec1}
The results are now presented in the form of theorems.\\

{\textbf {Theorem 2.1.}} For the standardized unit-Gompertz distribution with pdf (\ref{pdfs}), we have
 \begin{equation}\label{2.4}
\mu_{n;k}^{(r)}=\frac{{\rm e}^{k}}{(n-1)!}\sum_{i=0}^{n-1}(-1)^{n-i-1}\binom{n\!-\!1}{i}k^{n+\frac{r}{\theta}-i-1}\,\Gamma\left(\!i\!-\!\frac{r}{\theta}\!+\!1,k\!\right),~k,n\geq 1, r\in\mathbb{N},
\end{equation}

where $\Gamma(\cdot,\cdot)$ is the incomplete gamma function given by
 \begin{equation}\label{inc}
\Gamma(a,x)=\int_x^\infty t^{a-1} {\rm e}^{-t} {\rm d}t.
\end{equation}

\vspace{.20cm}
{\textbf {Proof.}
In view of (\ref{mean}), we can write
\begin{equation*}
\mu_{n;k}^{(r)}=\frac{k^n}{(n-1)!}\int_{0}^{1}x^{r}[-\ln {F}(x)]^{n-1}[{F}(x)]^{k-1}f(x){\rm d}x.
\end{equation*}
or, equivalently, from (1.4)
\begin{equation*}
\mu_{n;k}^{(r)}=\frac{\theta k^n}{(n-1)!}\int_{0}^{1}x^{r-\theta -1}[-\ln {F}(x)]^{n-1}[{F}(x)]^{k}{\rm d}x
\end{equation*}
\begin{equation}\label{2.6}
\qquad=\frac{\theta k^n\,{\rm e}^{k}}{(n-1)!}\sum_{i=0}^{n-1}(-1)^{n-i-1}\binom{n-1}{i}\int_{0}^{1}x^{r-\theta(i+1)-1}~{\rm e}^{-k\,x^{-\theta}}{\rm d}x.
\end{equation}

Setting $w = kx^{-\theta},$ we can write (\ref{2.6}) as
\begin{equation*}
\mu_{n;k}^{(r)}=\frac{{\rm e}^k}{(n-1)!}\sum_{i=0}^{n-1}(-1)^{n-i-1}~{k^{n+\frac{r}{\theta}-i-1}}\binom{n-1}{i}\int_{k}^{\infty}w^{i-\frac{r}{\theta}}\,{\rm e}^{-w}{\rm d}w,
\end{equation*}
and using the integral formula (\ref{inc}), it follows (\ref{2.4}).

\vspace{0.25cm}
{\textbf{Remark 2.1:(a)}} Setting $r=1$ in (\ref{2.4}), we obtain
\begin{equation}\label{2.7}
\mu_{n;k}\equiv\mu_{n;k}^{(1)}=\frac{{\rm e}^{k}}{(n-1)!}\sum_{i=0}^{n-1}(-1)^{n-i-1}\binom{n\!-\!1}{i}k^{n+\frac{1}{\theta}-i-1}\,\Gamma\left(\!i\!-\!\frac{1}{\theta}\!+\!1,k\!\right),
\end{equation}
which is the first moment (mean) of the $n$-th lower $k$-record value.

\vspace{0.25cm}
In addition, one can obtain the expression for the variance of the $n$-th lower $k$-record value by using the relation
\begin{equation}\label{2.8}
\sigma_{n;k}^{2}=\mu_{n;k}^{(2)}-\mu_{n;k}^{2},
\end{equation}
where $\mu_{n;k}^{(2)}$ follows from (\ref{2.4}) when $k = 2$.

\vspace{0.25cm}
We have computed the values of the means ($\mu_{n;k}$) and variances ($\sigma_{n;k}^{2}$) using (\ref{2.7}) and (\ref{2.8}), respectively, for $n=1,...,6$, $\theta=0.75, 1.5(1.0)4.5$ and $k=1,2,3$ and these values are presented in Tables 1 and 2, respectively. It can be seen from Table 1  that the means  are decreasing with respect to (w.r.t.) $n$.  Moreover, the means are increasing w.r.t. $\theta$. One can also observe from Table 1  that the means are increasing w.r.t. $k$.

\vspace{0.25cm}
{\bf(b)} Setting $k=n=1$ in (\ref{2.7}) and (\ref{2.8}), we obtain
$$\mu_{1;1}\equiv E[Z]={\rm e}~\Gamma\left(\!1\!-\!\frac{1}{\theta},1\right)$$ and $$\sigma_{1;1}^{2}\equiv Var[Z]=\mu_{1;1}^{(2)}-\mu_{1;1}^{2}
={\rm e}\left[\Gamma\left(\!1\!-\!\frac{2}{\theta},1\right)-{\rm e}\left\{\Gamma\left(\!1\!-\!\frac{2}{\theta},1\right)\right\}^2\right],$$
respectively, where $Z\sim UG(0,1,\theta)$.

\vspace{0.25cm}
One can easily extend the results for the UG distribution with pdf (\ref{pdfs}) by using a proper transformation (Mazucheli $\it{et \, al.}$, 2019).
	
\vspace{0.25cm}
{\bf(c)} By setting $k=1$ in (\ref{2.4}), we get
 \begin{equation}
\mu_{n;1}^{(r)}=\frac{{\rm e}}{(n-1)!}\sum_{i=0}^{n-1}(-1)^{n-i-1}\binom{n\!-\!1}{i}\Gamma\left(\!i\!-\!\frac{r}{\theta}\!+\!1,1\!\right),~n\geq 1, r\in\mathbb{N},
\end{equation}
which is an explicit expression for the $r$-th single moment of the $n$-th lower record value from the standardized UG distribution.

\vspace{0.25cm}
The following theorem provides a recurrence relation for the single moments of lower $k$-record values.

\vspace{.25cm}
{\textbf {Theorem 2.2.}} For $k\geq1$, $n>1$ and $r\in\mathbb{N}$, we have
\begin{equation}\label{2.10}
\mu_{n;k}^{(r+\theta)}=\frac{\theta\,k}{r}\Big(\mu_{n-1;k}^{(r)}-\mu_{n;k}^{(r)}\Big).
\end{equation}

\vspace{.25cm}
{\textbf {Proof.} From (\ref{mean}), we can write
\begin{equation*}
\mu_{n;k}^{(r)}=\frac{k^n}{(n-1)!}\int_{0}^{1} x^{r}\,[-\ln {F}(x)]^{n-1}[{F}(x)]^{k-1}\,f(x){\rm d}x.
\end{equation*}
Integrating by parts treating $[{F}(x)]^{k-1}\,f(x)$ for integration and the rest of the integrand for differentiation, we get
\begin{equation*}
\mu_{n;k}^{(r)}=\frac{k^{n-1}}{(n-2)!}\int_{0}^{1} x^{r}\,[-\ln {F}(x)]^{n-2}[{F}(x)]^{k-1}f(x)\,{\rm d}x
\end{equation*}
\begin{equation}\label{2.11}
\qquad\quad-\frac{r\,k^{n-1}}{(n-1)!}\int_{0}^{1} x^{r-1}\,[-\ln {F}(x)]^{n-1}[{F}(x)]^{k}{\rm d}x.
\end{equation}
Using relation (\ref{2.3}) in (\ref{2.11}) and rearranging the terms yields (\ref{2.10}). This completes the proof.

\vspace{0.25cm}
{\textbf{Remark 2.2.}} Setting $k=1$ in (\ref{2.10}), we obtain the recurrence relation for the single moments of ordinary lower record values from the standardized UG distribution.
\newpage
\begin{table}[H]
	{\bf Table 1:} Means of the lower $k$-record values.
	\setlength{\tabcolsep}{1.3em}
	\begin{center}
			\begin{tabular}{|c|c|c|c|c|c|}
				\cline{2-6}
				\multicolumn{1}{c|}{} &\multicolumn{5}{c|}{$k=1$}\\
				\hline
				$n$ &$\theta=0.75$ & $\theta=1.5$ & $\theta=2.5$ & $\theta=3.5$ & $\theta=4.5$\\
				\hline
   1  &   0.51742  & 0.69698  & 0.79852  & 0.84926  & 0.87962\\
   2  &   0.31010  & 0.53535  & 0.68059  & 0.75736  & 0.80453\\
   3  &   0.20703  & 0.43771  & 0.60344  & 0.69511  & 0.75268\\
   4  &   0.14937  & 0.37299  & 0.54870  & 0.64966  & 0.71421\\
   5  &   0.11400  & 0.32701  & 0.50751  & 0.61462  & 0.68415\\
   6  &   0.09067  & 0.29260  & 0.47515  & 0.58650  & 0.65976\\
				\hline
				\end{tabular}
			\end{center}
			\end{table}
\vspace{-1.1cm}
\begin{table}[H]
	\setlength{\tabcolsep}{1.3em}
	\begin{center}
			\begin{tabular}{|c|c|c|c|c|c|}
				\cline{2-6}
				\multicolumn{1}{c|}{} &\multicolumn{5}{c|}{$k=2$}\\
				\hline
				$n$ &$\theta=0.75$ & $\theta=1.5$ & $\theta=2.5$ & $\theta=3.5$ & $\theta=4.5$\\
				\hline
   1  &  0.65781   &      0.79920  	&    0.87079    &     0.90479    &     0.92464\\
   2  &  0.46511   &      0.66800   &	 0.78089    &     0.83669    &     0.86989\\
   3  &  0.34774   &      0.57653   &    0.71461    &     0.78527    &     0.82798\\
   4  &  0.27144   &      0.50939   &    0.66352    &     0.74476    &     0.79459\\
   5  &  0.21911   &      0.45806   &    0.62271    &     0.71182    &     0.76714\\
   6  &  0.18161   &      0.41752   &    0.58922    &     0.68432    &     0.74402\\
				\hline
				\end{tabular}
			\end{center}
			\end{table}
\vspace{-1.2cm}
\begin{table}[H]
	\setlength{\tabcolsep}{1.3em}
	\begin{center}
			\begin{tabular}{|c|c|c|c|c|c|}
				\cline{2-6}
				\multicolumn{1}{c|}{} &\multicolumn{5}{c|}{$k=3$}\\
				\hline
				$n$ &$\theta=0.75$ & $\theta=1.5$ & $\theta=2.5$ & $\theta=3.5$ & $\theta=4.5$\\
				\hline
   1  &  0.73194  & 0.84802    &     0.90375   &      0.92959    &     0.94450\\
   2  &  0.56020  &	0.73862    &     0.83100   &	  0.87522    &     0.90112\\
   3  &  0.44435  & 0.65651    &     0.77393   &      0.83175    &     0.86606\\
   4  &  0.36271  & 0.59273    &     0.72781   &      0.79601    &     0.83697\\
   5  &  0.30303  & 0.54178    &     0.68962   &      0.76595    &     0.81229\\
   6  &  0.25803  & 0.50011    &     0.65736   &      0.74022    &     0.79099\\
				\hline
				\end{tabular}
			\end{center}
			\end{table}

\subsection{Product Moments}\label{sec1}
The product moment of lower $k$-record values from the UG$(0,1,\theta)$ distribution are reported below.

\vspace{0.20cm}
{\textbf {Theorem 2.3.}} For $k\geq1$, $1\leq m< n$ and $r,s\in\mathbb{N}$ , we have
\begin{equation*}
\mu_{m,n;k}^{(r,s)}=\frac{\theta\,k^{n-1}\,{\rm e}^{k}}{(m-1)!\,(n-m-1)!}\sum_{i=0}^{m-1}\sum_{j=0}^{n-m-1}\binom{m-1}{i}\binom{n-m-1}{j}
\end{equation*}
\vspace{-0.5cm}
\begin{equation}\label{2.12}
\qquad\qquad\times\,\frac{(-1)^{n-m-1-j+i}}{[r-\theta(n-1-i-j)]}\bigg\{\frac{\Gamma(j\!-\!\frac{s}{\theta}\!+\!1,k)}{ k^{j-s/\theta}}-\frac{\Gamma(n\!-\!\frac{s+r}{\theta}\!-\!i,k)}{k^{n-i-1-{(s+r)}/\theta}}\bigg\},
\end{equation}
where $\Gamma(\cdot,\cdot)$ is the incomplete gamma function given in (\ref{inc}).\\

{\textbf {Proof.} From (\ref{pmean}), we have
\begin{equation}\label{2.13}
\mu_{m,n;k}^{(r,s)}=\frac{k^n}{(m-1)! \, (n-m-1)!}\int_{0}^{1} y^s\, [{F}(y)]^{k-1}f(y)\,{\rm I}(y)\,{\rm d}y,
\end{equation}
where
\begin{equation}\label{2.14}
{\rm I}(y)=\int_{y}^{1}x^r\,[-\ln {F}(x)]^{m-1}\frac{f(x)}{{F}(x)}\,[-\ln{F}(y)+\ln{F}(x)]^{n-m-1}\,{\rm d}x.
\end{equation}
Using (\ref{2.3}), (\ref{2.14}) can be rewritten as
\begin{equation*}
{\rm I}(y)=\int_{y}^{1}x^{r-\theta-1}\,[-\ln {F}(x)]^{m-1}[-\ln{F}(y)+\ln{F}(x)]^{n-m-1}\,{\rm d}x
\end{equation*}
\begin{equation*}
\qquad=\theta\int_{y}^{1}x^{r-\theta-1}\,(x^{-\theta}-1)^{m-1}\,(y^{-\theta}-x^{-\theta})^{n-m-1}\,{\rm d}x
\end{equation*}
\begin{equation*}
\qquad=\theta\sum_{i=0}^{m-1}\sum_{j=0}^{n-m-1}(-1)^{n-m+i-j-1}\binom{m-1}{i}\binom{n-m-1}{j}\,y^{-j\theta}\int_{y}^{1}x^{r-\theta(n-i-j-1)-1}\,{\rm d}x
\end{equation*}
\begin{equation*}
\qquad=\theta\sum_{i=0}^{m-1}\sum_{j=0}^{n-m-1}(-1)^{n-m+i-j-1}\binom{m-1}{i}\binom{n-m-1}{j}\,y^{-j\theta}\frac{[1-y^{r-\theta(n-i-j-1)}]}{[r-\theta(n-i-j-1)]}.
\end{equation*}

Inserting the resultant expression of ${\rm I}(y)$ in (\ref{2.13}), we get
\begin{equation*}
\mu_{m,n;k}^{(r,s)}=\frac{k^n\,\theta}{(m\!-\!1)!(n\!-\!m\!-\!1)!}\sum_{i=0}^{m-1}\sum_{j=0}^{n-m-1}\!\binom{\!m\!-\!1\!}{i}\binom{\!n\!-\!m\!-\!1\!}{j}
\frac{(-1)^{n-m+i-j-1}}{[r\!-\!\theta(n\!-\!i\!-\!j\!-\!1)]}
\end{equation*}
\begin{equation*}
\hspace{0.6in}\times\int_{0}^{1} y^{s-j\theta}(1-y^{r-\theta(n-i-j-1)})[{F}(y)]^{k-1}f(y){\rm d}y
\end{equation*}
\begin{equation*}
\qquad\quad=\frac{k^n{\rm e}^{k}\theta^2}{(m\!-\!1)!(n\!-\!m\!-\!1)!}\sum_{i=0}^{m-1}\sum_{j=0}^{n-m-1}\!\binom{\!m\!-\!1\!}{i}\binom{\!n\!-\!m\!-\!1\!}{j}
\frac{(-1)^{n-m+i-j-1}}{[r\!-\!\theta(n\!-\!i\!-\!j\!-\!1)]}
\end{equation*}
\begin{equation*}
\hspace{0.6in}\times\left[\int_{0}^{1}\!y^{s-\theta(j+1)-1}~{\rm e}^{-k y^{-\theta}}{\rm d}y
-\int_{0}^{1}\!y^{r+s-\theta(n-i)-1}~{\rm e}^{-k y^{-\theta}}{\rm d}y\right].
\end{equation*}
Setting $t=k\,y^{-\theta}$ and using integral formula (\ref{inc}), it follows (\ref{2.12}).

\vspace{0.25cm}
{\textbf{Remark 2.3:}} Setting $r=s=1$ in (\ref{2.12}), we obtain
\begin{equation*}
\mu_{m,n;k}\equiv\mu_{m,n;k}^{(1,1)}=\frac{\theta\,k^{n-1}\,{\rm e}^{k}}{(m-1)!\,(n-m-1)!}\sum_{i=0}^{m-1}\sum_{j=0}^{n-m-1}\binom{m-1}{i}\binom{n-m-1}{j}
\end{equation*}
\begin{equation*}
\qquad\qquad\qquad\qquad\times\,\frac{(-1)^{n-m-1-j+i}}{[1-\theta(n-1-i-j)]}\bigg\{\frac{\Gamma(j\!-\!\frac{1}{\theta}\!+\!1,k)}{ k^{j-1/\theta}}-\frac{\Gamma(n\!-\!\frac{2}{\theta}\!-\!i,k)}{k^{n-i-1-{2}/\theta}}\bigg\},
\end{equation*}
\begin{equation}
\hspace{4in}k\geq1,~1\leq m< n,
\end{equation}
which is the $(1,1)$-th (say simple) product moment of the $Z_{m(k)}$ and $Z_{n(k)}$. The simple product moments are used for evaluating the covariances, i.e.
\begin{equation*}
\sigma_{m,n;k}\equiv Cov(Z_{m(k)},Z_{n(k)})=\mu_{m,n;k}-\mu_{m;k}~\mu_{n;k}.
\end{equation*}

We have computed the values of covariances $(\sigma_{m,n;k})$ for different choices of $m, n, k$ and $\theta$, and the results are presented in Table 2.

\vspace{0.25cm}
The following theorem provides a recurrence relation for the product moments of $k$-record values from the standardized UG distribution with pdf (\ref{pdfs}).
	
\vspace{0.25cm}
{\textbf {Theorem 2.4.}} For the standardized unit-Gompertz distribution with the pdf given in (\ref{pdfs}), $k\geq 1$, $n>m>1$ and $r,s\in \mathbb{N}$,
we have
\begin{equation}\label{t2}
\mu_{m,n;k}^{(r,s+\theta)}=\frac{\theta\,k}{s}\Big(\mu_{m,n-1;k}^{(r,s)}-\mu_{m,n;k}^{(r,s)}\Big).
\end{equation}

\vspace{0.25cm}
{\textbf {Proof.} From (\ref{pmean}), we have
\begin{equation}
\mu_{m,n;k}^{(r,s)}=\frac{k^n}{(m-1)! \, (n-m-1)!}\int_{0}^{1} x^r\, [-\ln {F}(x)]^{m-1}\frac{f(x)}{{F}(x)}{\rm I}(x)\,{\rm d}x,
\end{equation}
where
\begin{equation*}
{\rm I}(x)=\int_{0}^{x}y^s\,[\ln{F}(x)-\ln{F}(y)]^{n-m-1}\,[{F}(y)]^{k-1}f(y)\,{\rm d}y.
\end{equation*}
Integrating by parts treating $[{F}(y)]^{k-1}\,f(y)$ for integration and the rest of the integrand for differentiation, we get
\begin{equation*}
{\rm I}(x)=\frac{n-m-1}{k}\int_{0}^{x}y^s\,[\ln{F}(x)-\ln{F}(y)]^{n-m-2}\,[{F}(y)]^{k-1}f(y)\,{\rm d}y
\end{equation*}
\begin{equation}\label{2.18}
\hspace{1.2cm}-\frac{s}{k}\int_{0}^{x}y^{s-1}\,[\ln{F}(x)-\ln{F}(y)]^{n-m-1}[{F}(y)]^{k}\,{\rm d}y.
\end{equation}
Now, using relation (\ref{2.3}), in (\ref{2.18}) and rearranging the terms, leads to (\ref{t2}).

\vspace{0.25cm}
{\textbf{Remark 2.4.}}~{\bf(a)} Setting $r=0$ in (\ref{t2}), the recurrence relation for the single moments of $k$-record values  from the standardized UG distribution given in (\ref{2.10}) will be recovered.

\vspace{0.25cm}
{\bf(b)} At $k = 1$ in (\ref{t2}), we found the recurrence relation for the product moments of the lower $k$-record values from the standardized UG distribution.

\vspace{0.25cm}
{\textbf{Remark 2.5:}} The results of this section can easily be extended for the unstandardized UG distribution with pdf (\ref{pdf}) by applying proper transformation (see Kumar $\it{et \, al.}$, 2020).
			
\begin{table}[H]
{\bf Table 2:} Variances and covariances of the lower $k$-record values.\vspace{-0.1cm}
\setlength{\tabcolsep}{1.3em}
			\begin{center}
				\begin{tabular}{|c|c|c|c|c|c|c|}
				\cline{3-7}
				\multicolumn{2}{c|}{} &\multicolumn{5}{c|}{$k=1$}\\
				\hline
				$m$ &$n$ &$\theta=0.75$ & $\theta=1.5$ & $\theta=2.5$ & $\theta=3.5$ & $\theta=4.5$\\
				\hline
   1  &   1 & 0.05955 & 0.03164 & 0.01604 & 0.00955 & 0.00631\\
   1  &   2 & 0.02824 & 0.01851 & 0.01023 & 0.00633 & 0.00427\\
   1  &   3 & 0.01518 & 0.01196 & 0.00712 & 0.00454 & 0.00312\\
   1  &   4 & 0.00901 & 0.00833 & 0.00528 & 0.00346 & 0.00241\\
   1  &   5 & 0.00579 & 0.00614 & 0.00410 & 0.00275 & 0.00194\\
   1  &   6 & 0.00395 & 0.00472 & 0.00330 & 0.00226 & 0.00161\\
   2  &   2 & 0.03108 & 0.02350 & 0.01386 & 0.00882 & 0.00605\\
   2  &   3 & 0.01645 & 0.01507 & 0.00959 & 0.00630 & 0.00440\\
   2  &   4 & 0.00964 & 0.01043 & 0.00708 & 0.00478 & 0.00339\\
   2  &   5 & 0.00613 & 0.00765 & 0.00548 & 0.00379 & 0.00272\\
   2  &   6 & 0.00415 & 0.00586 & 0.00440 & 0.00310 & 0.00225\\
   3  &   3 & 0.01474 & 0.01543 & 0.01041 & 0.00701 & 0.00497\\
   3  &   4 & 0.00856 & 0.01063 & 0.00766 & 0.00531 & 0.00381\\
   3  &   5 & 0.00540 & 0.00777 & 0.00591 & 0.00420 & 0.00306\\
   3  &   6 & 0.00363 & 0.00593 & 0.00474 & 0.00343 & 0.00253\\
   4  &   4 & 0.00730 & 0.01025 & 0.00777 & 0.00551 & 0.00401\\
   4  &   5 & 0.00458 & 0.00747 & 0.00599 & 0.00435 & 0.00321\\
   4  &   6 & 0.00306 & 0.00569 & 0.00480 & 0.00355 & 0.00265\\
   5  &   5 & 0.00387 & 0.00706 & 0.00593 & 0.00439 & 0.00327\\
   5  &   6 & 0.00258 & 0.00538 & 0.00474 & 0.00358 & 0.00270\\
   6  &   6 & 0.00220 & 0.00506 & 0.00464 & 0.00356 & 0.00271\\
				\hline
				\end{tabular}
			\end{center}
			\end{table}

\begin{table}[H]
{\bf Table 2:} continued...\vspace{-0.1cm}
\setlength{\tabcolsep}{1.3em}
			\begin{center}
				\begin{tabular}{|c|c|c|c|c|c|c|}
				\cline{3-7}
				\multicolumn{2}{c|}{} &\multicolumn{5}{c|}{$k=2$}\\
				\hline
				$m$ &$n$ &$\theta=0.75$ & $\theta=1.5$ & $\theta=2.5$ & $\theta=3.5$ & $\theta=4.5$\\
				\hline
   1  &   1 & 0.04440 & 0.01909 & 0.00873 & 0.00496 & 0.00318\\
   1  &   2 & 0.02649 & 0.01313 & 0.00637 & 0.00371 & 0.00242\\
   1  &   3 & 0.01691 & 0.00953 & 0.00487 & 0.00290 & 0.00191\\
   1  &   4 & 0.01143 & 0.00722 & 0.00386 & 0.00234 & 0.00156\\
   1  &   5 & 0.00808 & 0.00566 & 0.00315 & 0.00195 & 0.00131\\
   1  &   6 & 0.00594 & 0.00456 & 0.00263 & 0.00165 & 0.00112\\
   2  &   2 & 0.03425 & 0.01888 & 0.00959 & 0.00570 & 0.00376\\
   2  &   3 & 0.02171 & 0.01365 & 0.00731 & 0.00445 & 0.00297\\
   2  &   4 & 0.01458 & 0.01031 & 0.00579 & 0.00359 & 0.00242\\
   2  &   5 & 0.01027 & 0.00806 & 0.00471 & 0.00298 & 0.00203\\
   2  &   6 & 0.00751 & 0.00648 & 0.00393 & 0.00252 & 0.00173\\
   3  &   3 & 0.02210 & 0.01535 & 0.00858 & 0.00532 & 0.00359\\
   3  &   4 & 0.01476 & 0.01156 & 0.00678 & 0.00429 & 0.00292\\
   3  &   5 & 0.01035 & 0.00902 & 0.00551 & 0.00355 & 0.00245\\
   3  &   6 & 0.00755 & 0.00724 & 0.00459 & 0.00300 & 0.00209\\
   4  &   4 & 0.01392 & 0.01195 & 0.00729 & 0.00469 & 0.00323\\
   4  &   5 & 0.00972 & 0.00931 & 0.00592 & 0.00388 & 0.00270\\
   4  &   6 & 0.00707 & 0.00746 & 0.00493 & 0.00328 & 0.00230\\
   5  &   5 & 0.00890 & 0.00929 & 0.00612 & 0.00407 & 0.00285\\
   5  &   6 & 0.00646 & 0.00744 & 0.00509 & 0.00344 & 0.00243\\
   6  &   6 & 0.00585 & 0.00729 & 0.00515 & 0.00353 & 0.00252\\
				\hline
				\end{tabular}
			\end{center}
			\end{table}

\begin{table}[H]
{\bf Table 2:} continued...\vspace{-0.1cm}
\setlength{\tabcolsep}{1.3em}
			\begin{center}
				\begin{tabular}{|c|c|c|c|c|c|c|}
				\cline{3-7}
				\multicolumn{2}{c|}{} &\multicolumn{5}{c|}{$k=3$}\\
				\hline
				$m$ &$n$ &$\theta=0.75$ & $\theta=1.5$ & $\theta=2.5$ & $\theta=3.5$ & $\theta=4.5$\\
				\hline
   1  &   1 & 0.03320 & 0.01281 & 0.00557 & 0.00309 & 0.00196\\
   1  &   2 & 0.02213 & 0.00954 & 0.00434 & 0.00246 & 0.00157\\
   1  &   3 & 0.01544 & 0.00736 & 0.00349 & 0.00201 & 0.00130\\
   1  &   4 & 0.01119 & 0.00585 & 0.00288 & 0.00168 & 0.00110\\
   1  &   5 & 0.00837 & 0.00476 & 0.00242 & 0.00144 & 0.00095\\
   1  &   6 & 0.00644 & 0.00396 & 0.00207 & 0.00125 & 0.00083\\
   2  &   2 & 0.03112 & 0.01464 & 0.00692 & 0.00398 & 0.00258\\
   2  &   3 & 0.02162 & 0.01127 & 0.00555 & 0.00325 & 0.00213\\
   2  &   4 & 0.01561 & 0.00894 & 0.00457 & 0.00272 & 0.00180\\
   2  &   5 & 0.01164 & 0.00726 & 0.00384 & 0.00232 & 0.00154\\
   2  &   6 & 0.00893 & 0.00603 & 0.00329 & 0.00201 & 0.00135\\
   3  &   3 & 0.02357 & 0.01334 & 0.00681 & 0.00405 & 0.00267\\
   3  &   4 & 0.01697 & 0.01056 & 0.00560 & 0.00338 & 0.00225\\
   3  &   5 & 0.01262 & 0.00857 & 0.00470 & 0.00289 & 0.00194\\
   3  &   6 & 0.00966 & 0.00710 & 0.00402 & 0.00250 & 0.00169\\
   4  &   4 & 0.01693 & 0.01138 & 0.00623 & 0.00382 & 0.00256\\
   4  &   5 & 0.01257 & 0.00923 & 0.00523 & 0.00325 & 0.00220\\
   4  &   6 & 0.00960 & 0.00764 & 0.00447 & 0.00282 & 0.00192\\
   5  &   5 & 0.01206 & 0.00951 & 0.00555 & 0.00350 & 0.00238\\
   5  &   6 & 0.00920 & 0.00787 & 0.00474 & 0.00303 & 0.00208\\
   6  &   6 & 0.00866 & 0.00792 & 0.00491 & 0.00317 & 0.00219\\
				\hline
				\end{tabular}
			\end{center}
			\end{table}

\section{Best Linear Unbiased Estimation}\label{sec1}
In this section, we obtain the BLUEs of location and scale parameters of the three parameter UG distribution. Let $Z_{1(k)}, Z_{2(k)}, \ldots, Z_{n(k)}$ be the first $n$ lower \emph{k}-record values arising from the standardized UG distribution with pdf given in (\ref{pdfs}) and let $R_{i(k)}=\sigma Z_{i(k)}+\mu$ for $i= 1,2,\ldots,n$.
Then $R_{1(k)}, R_{2(k)}, \ldots, R_{n(k)}$ are
 the corresponding lower \emph{k}-record values from the UG  distribution with location parameter $\mu$ and scale parameter $\sigma$. If we denote $E\left(Z_{n(k)}\right)=\alpha_{n(k)}$, $Var\left(Z_{n(k)}\right)=\beta_{n,n(k)}$ and $Cov\left(Z_{m(k)},Z_{n(k)}\right)=\beta_{m,n(k)}$, then we have for $ n \geq 1$
\[E\left(R_{n(k)}\right)=\mu+\sigma\alpha_{n(k)}~,\]
\[Var\left(R_{n(k)}\right)=\sigma^{2}\beta_{n,n(k)}\]
and for  $m<n$
\[Cov\left(R_{m(k)},R_{n(k)}\right)=\sigma^{2}\beta_{m,n(k)}.\]
Suppose $\textbf{R}_{n(k)}=\left(R_{1(k)},R_{2(k)},...,R_{n(k)}\right)^{T}$ denotes the vector of lower \emph{k}-record values. Then
\[E\left(\textbf{R}_{n(k)}\right)=\mu\textbf{1}+\sigma \bm{\alpha},\]
where $\bm{\alpha}=\left(\alpha_{1(k)},\alpha_{2 (k)},...,\alpha_{n(k)}\right)^{T}$ and \textbf{1} is a column vector of $n$ ones. The variance-covariance matrix of $\textbf{R}_{n(k)}$ is given by
\[ D\left(\textbf{R}_{n(k)}\right)=\textbf{B}\sigma^{2},\]
where $\textbf{B}=\left((\beta_{i,j(k)},~1\leq i,j\leq n)\right)$. Following the generalized least-squares approach, the BLUEs of $\mu$ and $\sigma$ are given, respectively, by (see e.g. Balakrishnan and Cohen, 1991)
\begin{eqnarray}\label{3.3}
\bm{\mu^{\ast}}&=&\frac{\bm{\alpha}^{T}\textbf{B}^{-1}\bm{\alpha\textbf{1}^{T}}\textbf{B}^{-1}-\bm{\alpha}^{T}\textbf{B}^{-1}\textbf{1}\bm{\alpha^{T}}\textbf{B}^{-1}}{(\bm{\alpha}^{T}\textbf{B}^{-1}\bm{\alpha})(\textbf{1}^{T}\textbf{B}^{-1}\textbf{1})-(\bm{\alpha}^{T}\textbf{B}^{-1}\textbf{1})^{2}}\textbf{R}_{n(k)}\nonumber\\
&=& \sum_{i=1}^{n}a_{i}R_{i(k)}
\end{eqnarray}
and
\begin{eqnarray}\label{3.4}
\bm{\sigma^{\ast}}&=&\frac{\textbf{1}^{T}\textbf{B}^{-1}\textbf{1}\bm{\alpha}^{T}\textbf{B}^{-1}-\textbf{1}^{T}\textbf{B}^{-1}\bm{\alpha}\textbf{1}^{T}\textbf{B}^{-1}}{(\bm{\alpha}^{T}\textbf{B}^{-1}\bm{\alpha})(\textbf{1}^{T}\textbf{B}^{-1}\textbf{1})-(\bm{\alpha}^{T}\textbf{B}^{-1}\textbf{1})^{2}}\textbf{R}_{n(k)}\nonumber\\
&=& \sum_{i=1}^{n}b_{i}R_{i(k)}.
\end{eqnarray}
Furthermore, the variances and covariance of the above estimators are given by
\begin{equation}\label{3.5}
Var(\mu^{\ast})=\sigma^{2}\left(\frac{\bm{\alpha}^{T}\textbf{B}^{-1}\bm{\alpha}}{(\bm{\alpha}^{T}\textbf{B}^{-1}\bm{\alpha})(\textbf{1}^{T}\textbf{B}^{-1}\textbf{1})-(\bm{\alpha}^{T}\textbf{B}^{-1}\textbf{1})^{2}}\right)=\sigma^2 V_1,
\end{equation}
\begin{equation}\label{3.6}
Var(\sigma^{\ast})=\sigma^{2}\left(\frac{\textbf{1}^{T}\textbf{B}^{-1}\textbf{1}}{(\bm{\alpha}^{T}\textbf{B}^{-1}\bm{\alpha})(\textbf{1}^{T}\textbf{B}^{-1}\textbf{1})-(\bm{\alpha}^{T}\textbf{B}^{-1}\textbf{1})^{2}}\right)=\sigma^2 V_2    \end{equation}
and
\begin{equation}\label{3.6a}
Cov(\mu^{\ast},\sigma^{\ast})=\sigma^{2}\left(\frac{-\bm{\alpha}^{T}\textbf{B}^{-1}\textbf{1}}{(\bm{\alpha}^{T}\textbf{B}^{-1}\bm{\alpha})(\textbf{1}^{T}\textbf{B}^{-1}\textbf{1})-(\bm{\alpha}^{T}\textbf{B}^{-1}\textbf{1})^{2}}\right)=\sigma^2V_3.
\end{equation}
\par
By making use of the values of means, variances and covariances presented in Tables 1 and 2, we have calculated the coefficients $a_{i}$ and $b_{i},~i=1,2,\ldots,n$ of BLUEs of $\mu$ and $\sigma$ using (\ref{3.3}) and (\ref{3.4}) for $n=2(1)6;~\theta=0.75, 1.5(1.0)4.5$ and $k=1,2,3$  and the results are given in Tables 3 and 4. The values of $V_1=\dfrac{Var(\mu^{\ast})}{\sigma^{2}},~~V_2=\dfrac{Var(\sigma^{\ast})}{\sigma^{2}} ~\mbox{and} ~~V_3=\dfrac{Cov(\mu^{\ast},\sigma^{\ast})}{\sigma^{2}}$ for $n=2(1)6;~\theta=0.75, 1.5(1.0)4.5$ and $k=1,2,3$  are also computed from equations (\ref{3.5}), (\ref{3.6}) and (\ref{3.6a}) and are incorporated in Table 5.

\vspace{0.25cm}
Based on the  BLUEs of the location and scale  parameters, the
confidence intervals (CIs) for $\mu$ and $\sigma$ can be constructed
through the pivotal quantities given by
\begin{equation*}\label{quant}
T_1=\frac{\mu^*-\mu}{\sigma^*\sqrt{V_1}}
\hspace{1cm}\textrm{and}\hspace{1cm} T_2=\frac{\sigma^*-\sigma}{\sigma\sqrt{V_2}}.
\end{equation*}

Constructing such CIs requires then the percentage points of $T_1$
and $T_2$ which  can be computed by using the BLUEs
$\mu^*$ and $\sigma^*$ via Monte
Carlo method. In Table 6, we have determined the percentage points
of $T_1$ and $T_2$  based on 10000 runs and different values of $n$
and $\theta$. Based on these simulated percentage points, we can
determine a $100(1-\gamma)\%$ CI for $\mu$  through the pivotal
quantity $T_1$ as follows
\begin{equation}\label{cimu}
P\bigg{(}\mu^*-\sigma^*
\sqrt{V_{1}}T_1(1-\gamma/2)\leq \mu \leq
\mu^*-\sigma^*
\sqrt{V_{1}}T_1(\gamma/2)\bigg{)}=1-\gamma,
\end{equation}
where $T_1(\tau)$ is the left percentage point of $T_1$ at $\tau$, i.e. $P(T_1<T_1(\tau))=\tau$.

Similarly, a $100(1-\gamma)\%$  CI for $\sigma$ can be constructed
through the pivotal quantity $T_2$ as follows
\begin{equation}\label{cisigma}
P\left(\frac{\sigma^*}{1+\sqrt{V_{2}}T_2(1-\gamma/2)}\leq
\sigma \leq
\frac{\sigma^*}{1+\sqrt{V_{2}}T_2(\gamma/2)}\right)=1-\gamma,
\end{equation}
where $T_2(\tau)$ is the left percentage point of $T_2$ at $\tau$, i.e. $P(T_2<T_2(\tau))=\tau$.

\section{Best Linear Invariant Estimation}\label{sec1}

Now let us consider the best linear invariant estimators of $\mu$ and $\sigma$. Based on the results of Mann (1969), the BLIEs of $\mu$ and $\sigma$ are given by
 (see also Arnold et al., 1998, p. 143)
\[\widetilde{\mu}=\mu^{\ast}-\frac{V_3}{1+V_2}\sigma^{\ast} \quad \textrm{and} \quad \widetilde{\sigma}=\frac{\sigma^{\ast}}{1+V_2},\]
where $V_1=\frac{1}{\sigma^2}Var(\mu^{\ast}),$
$V_2=\frac{1}{\sigma^2}Var(\sigma^{\ast})$ and
$$V_3=\frac{1}{\sigma^2}Cov(\mu^{\ast},
\sigma^{\ast}).$$
Furthermore, the variances of these
BLIEs are given by (see Arnold et al., 1998, p. 143)
\[Var(\widetilde{\mu})=\sigma^2\left(V_1-\frac{V_3^2(2+V_2)}{(1+V_2)^2}\right) \quad \textrm{and} \quad Var(\widetilde{\sigma})=\frac{\sigma^2V_2}{(1+V_2)^2}.\]
Moreover, we have
\[
Cov(\widetilde{\mu},\widetilde{\sigma})=\frac{\sigma^2 V_3}{(1+V_2)^2}.
\]

Based on the  BLIEs, we can again construct CIs for the location and
scale parameters through pivotal quantities given by
\begin{equation*}\label{quanti}
T_3=\frac{\widetilde{\mu}-\mu}{\widetilde{\sigma}\sqrt{V_1-\frac{V_3^2(2+V_2)}{(1+V_2)^2}}}
\hspace{1cm}\textrm{and}\hspace{1cm} T_4=\frac{\widetilde{\sigma}-\sigma}{\sigma\frac{\sqrt{V_2}}{1+V_2}}.
\end{equation*}
Table 7 presents the percentage points of $T_3$ and $T_4$ based on
10000 runs and different choices of $n$ and $\theta$. With the BLIEs
and the use of Table 7, we can determine a $100(1-\gamma)\%$  CI for
$\mu$  through the pivotal quantity  $T_3$ as
\begin{eqnarray}\nonumber\label{cimu2}
P\bigg{(}\widetilde{\mu}-\widetilde{\sigma} \sqrt{V_1-\frac{V_3^2(2+V_2)}{(1+V_2)^2}}T_3(1-\gamma/2)\leq \mu \leq
\widetilde{\mu}-\widetilde{\sigma} \sqrt{V_1-\frac{V_3^2(2+V_2)}{(1+V_2)^2}}T_3(\gamma/2)\bigg{)}=1-\gamma,\\
\end{eqnarray}
where $T_3(\tau)$ is the left percentage point of $T_3$ at $\tau$, i.e. $P(T_3<T_3(\tau))=\tau$.

\vspace{0.25cm}
Similarly, we can determine a $100(1-\gamma)\%$  CI for $\sigma$,
through the pivotal quantity  $T_4$ as
\begin{equation}\label{cisigma2}
P\left(\frac{\widetilde{\sigma}}{1+\frac{\sqrt{V_2}}{1+V_2}T_4(1-\gamma/2)}\leq
\sigma \leq
\frac{\widetilde{\sigma}}{1+\frac{\sqrt{V_2}}{1+V_2}T_4(\gamma/2)}\right)=1-\gamma,
\end{equation}
where $T_4(\tau)$ is the left percentage point of $T_4$ at $\tau$, i.e. $P(T_4<T_4(\tau))=\tau$.

\vspace{0.25cm}
Now, let us compare the  BLUEs and BLIEs using the relative
efficiency criterion (REC). Since the mean squared errors (MSEs) of
BLUEs are equal to their corresponding variances, we have
\[{\rm MSE}(\mu^{\ast})=\sigma^2 V_1
\quad \textrm{and} \quad {\rm MSE}(\sigma^{\ast})=\sigma^2
V_2.\] On the other hand, the  MSEs of BLIEs of $\mu$ and $\sigma$
can be obtained as
\[{\rm MSE}(\widetilde{\mu})=\sigma^2\left(V_1-\frac{V_3^2}{1+V_2}\right)
\quad \textrm{and} \quad
{\rm MSE}(\widetilde{\sigma})=\frac{\sigma^2V_2}{1+V_2}.\]
Therefore, we can readily obtain the RECs of the BLIEs of $\mu$ and
$\sigma$ w.r.t. their corresponding BLUEs as follows
\begin{eqnarray}\label{rec1}
REC(\widetilde{\mu},
\mu^{\ast})&=&\dfrac{{\rm MSE}(\mu^{\ast})}{{\rm MSE}(\widetilde{\mu})}=
\dfrac{V_1}{V_1-\dfrac{V_3^2}{1+V_2}}\geq 1,\\ \label{rec2}
REC(\widetilde{\sigma},
\sigma^{\ast})&=&\dfrac{{\rm MSE}(\sigma^{\ast})}{{\rm MSE}(\widetilde{\sigma})}=
1+V_2\geq 1.
\end{eqnarray}
Therefore, both of the BLIEs of $\mu$ and $\sigma$ perform better
than the corresponding BLUEs in terms of MSE.
\begin{table}[H]
{\bf Table 3:} Coefficients for the BLUEs of $\mu$.\vspace{-0.1cm}
\begin{center}
\begin{tabular}{|c|c|c|c|c|c|c|c|c|}
\hline
\multicolumn{1}{|c}  {$k$} &\multicolumn{1}{|c|}{$\theta$} & $n$  & \multicolumn{6}{c|}{$a_i$,~$i =1,2,\ldots,n$} \\
\hline
 1 & 0.75 & 2 & $-1.49572$ &  $~~2.49572$  & & & &\\
   &      & 3 & $-0.54644$ &  $-0.36301$ & ~~1.90945 & & &\\
   &      & 4 & $-0.28907$ &  $-0.15222$ & $-0.32115$ & ~~1.76243 & &\\
   &      & 5 & $-0.18275$ &  $-0.07621$ & $-0.15712$ & $-0.30267$ & ~~1.71876  &\\
   &      & 6 & $-0.12822$ &  $-0.04287$ & $-0.08656$ & $-0.16479$ & $-0.29167$ &  1.71411\\
\cline{2-9}
     & $1.5$   & $2$   & $-3.31207$    & $~~4.31207$ & $$      & $$       & $$          & $$ \\
       &        & $3$ & $-1.47021$  & $-0.57910$ & $~~3.04931$   & $$          & $$          & $$ \\
\hline
		\end{tabular}
	\end{center}
\end{table}

\begin{table}[H]
{\bf Table 3:} continued...\vspace{-0.1cm}
\begin{center}
\begin{tabular}{|c|c|c|c|c|c|c|c|c|}
\hline
\multicolumn{1}{|c}  {$k$} &\multicolumn{1}{|c|}{$\theta$} & $n$  & \multicolumn{6}{c|}{$a_i$,~$i =1,2,\ldots,n$} \\
\hline
$2$    & $0.75$ & $2$ &   $-2.41364$ & ~~3.41364  &  & & &\\
       &        & $3$ &   $-0.96028$ & $-0.42581$ & ~~2.38609  &  & &\\
	   &        & $4$ &   $-0.53437$ & $-0.20033$ & $-0.34299$ & ~~2.07769 & &\\
	   &        & $5$ &   $-0.34803$ & $-0.10945$ & $-0.18506$ & $-0.29999$ & ~~1.94253 &\\
	   &        & $6$ &   $-0.24914$ & $-0.06574$ & $-0.10983$ & $-0.17667$ & $-0.27243$ & 1.87381\\
\cline{2-9}
     & $1.5$  & $2$ &   $-5.09159$ & ~~6.09159  &  &  &  &\\
       &        & $3$ &   $-2.32415$ & $-0.64523$ & ~~3.96938  &  &  &\\
\hline
		\end{tabular}
	\end{center}
\end{table}

\begin{table}[H]
{\bf Table 3:} continued...\vspace{-0.1cm}
\begin{center}
\begin{tabular}{|c|c|c|c|c|c|c|c|c|}
\hline
\multicolumn{1}{|c}  {$k$} &\multicolumn{1}{|c|}{$\theta$} & $n$  & \multicolumn{6}{c|}{$a_i$,~$i =1,2,\ldots,n$} \\
\hline
$3$    & $0.75$ & $2$ &   $-3.26186$ & ~~4.26186  &  &  &  &\\
       &        & $3$ &   $-1.35620$ & $-0.46879$ & ~~2.82499 &  &  &\\
	   &        & $4$ &   $-0.77708$ & $-0.23454$ & $-0.36104$ & ~~2.37266  &  &\\
	   &        & $5$ &   $-0.51549$ & $-0.13477$ & $-0.20583$ & $-0.30486$ & ~~2.16094 &\\
	   &        & $6$ &   $-0.37334$ & $-0.08432$ & $-0.12777$ & $-0.18823$ & $-0.26917$ &  2.04283\\
\cline{2-9}
     & $1.5$  & $2$ &   $-6.75175$ & ~~7.75175  &  &  &  &\\
       &        & $3$ &   $-3.12849$ & $-0.69878$ & ~~4.82727 & & &\\
\hline
		\end{tabular}
	\end{center}
\end{table}

\begin{table}[H]
{\bf Table 4:} Coefficients for the BLUEs of $\sigma$.\vspace{-0.1cm}
\begin{center}
\begin{tabular}{|c|c|c|c|c|c|c|c|c|}
\hline
\multicolumn{1}{|c}  {$k$} &\multicolumn{1}{|c|}{$\theta$} & $n$  & \multicolumn{6}{c|}{$b_i$,~$i =1,2,\ldots,n$} \\
\hline
$1$    & $0.75$ & $2$ &    4.82335 & $-4.82335$ & & & &\\
       &        & $3$ &    2.88363 & ~~1.01806 & $-3.90169$ & & &\\
	   &        & $4$ &    2.30214 &  ~~0.54180 & ~~1.13797 & $-3.98191$ & &\\
	   &        & $5$ &    2.02635 & ~~0.34464 & ~~0.71248 & ~~1.37491 & $-4.45838$ &\\
	   &        & $6$ &    1.85947 & ~~0.24260 & ~~0.49651 & ~~0.95292  & ~~1.69456 & $-5.24606$\\
\cline{2-9}
     & $1.5$  & $2$ &    6.18678 & $-6.18678$ & & &\\
       &        & $3$ &    3.47607 & ~~1.01169 & $-4.48776$& & &\\
\hline
		\end{tabular}
	\end{center}
\end{table}

\begin{table}[H]
{\bf Table 4:} continued...\vspace{-0.1cm}
\begin{center}
\begin{tabular}{|c|c|c|c|c|c|c|c|c|}
\hline
\multicolumn{1}{|c}  {$k$} &\multicolumn{1}{|c|}{$\theta$} & $n$  & \multicolumn{6}{c|}{$b_i$,~$i =1,2,\ldots,n$} \\
\hline
$2$    & $0.75$ & $2$ &    5.18941 & $-5.18941$ &  & &  & \\
       &        & $3$ &    2.93243 &  ~~0.77302&  $-3.70545$ & & & \\
	   &        & $4$ &    2.24503 &  ~~0.40912 &  ~~0.69908 &  $-3.35324$ &  & \\
	   &        & $5$ &    1.92760 &  ~~0.25430 &  ~~0.43004 &  ~~0.69722 & $-3.30916$&\\
	   &        & $6$ &    1.74729 &  ~~0.17460 &  ~~0.29286 &  ~~0.47235 &  ~~0.72967 &  $-3.41676$\\
\cline{2-9}
     & $1.5$  & $2$ &    7.62211 & $-7.62211$ & & & &\\
       &        & $3$ &    4.12677 & ~~0.88667 & $-5.01344$ & & &\\
\hline
		\end{tabular}
	\end{center}
\end{table}

\begin{table}[H]
{\bf Table 4:} continued...\vspace{-0.1cm}
\begin{center}
\begin{tabular}{|c|c|c|c|c|c|c|c|c|}
\hline
\multicolumn{1}{|c}  {$k$} &\multicolumn{1}{|c|}{$\theta$} & $n$  & \multicolumn{6}{c|}{$b_i$,~$i =1,2,\ldots,n$} \\
\hline
$3$    & $0.75$ & $2$ &    5.82269 &  $-5.82269$& & & &\\
       &        & $3$ &    3.19001 & ~~0.71272  &$-3.90273$ & & &\\
	   &        & $4$ &    2.37369 & ~~0.38252  &~~0.58824  &$-3.34446$ & &\\
	   &        & $5$ &    1.99440 & ~~0.23786  &~~0.36319  &~~0.53785  & $-3.13329$ &\\
	   &        & $6$ &    1.78077 & ~~0.16203  &~~0.24587  &~~0.36258  &~~0.51884   &$-3.07010$\\
\cline{2-9}
     & $1.5$  & $2$ &    9.14102 & $-9.14102$ & & & &\\
       &        & $3$ &    4.84800 & ~~0.87160  &$-5.71959$& & &\\
\hline
		\end{tabular}
	\end{center}
\end{table}
			
\begin{table}[H]
{\bf Table 5:} Variances and covariance of the BLUEs of $\mu$ and $\sigma$ in terms of $\sigma^2$ and $V_4$.\vspace{-0.1cm}
\begin{center}
\begin{tabular}{|c|c|c|c|c|c|c|}
				\hline
$k$   & $\theta$ & $n$ &  $V_1$   &  $V_2$  & $V_3$ & $V_4$\\
				\hline
$1$     & $0.75$ & $2$ &  0.11604 &0.79487 &$-0.26023$	&	$-0.08954$\\
		  &        & $3$ &  0.03236 &0.44551 &$-0.08925$	&	$-0.02513$\\
	      &        & $4$ &  0.01277 &0.34546 &$-0.04497$	&	$-0.01015$\\
		  &        & $5$ &  0.00604 &0.30019 &$-0.02752$	&	$-0.00501$\\
		  &        & $6$ &  0.00323 &0.27389 &$-0.01892$	&	$-0.00284$\\
				\cline{2-7}
        & $1.5$  & $2$ &  0.25555 &0.69414 &$-0.40253$	&	$-0.07950$\\
		  &        & $3$ &  0.09083 &0.33736 &$-0.16011$	&	$-0.02608$\\
	      &        & $4$ &  0.04557 &0.23064 &$-0.09061$	&	$-0.01208$\\
		  &        & $5$ &  0.02699 &0.18181 &$-0.06049$	&	$-0.00664$\\
		  &        & $6$ &  0.01772 &0.15422 &$-0.04449$	&	$-0.00409$\\
				\cline{2-7}
        & $2.5$  & $2$ &  0.37000 &0.67848 &$-0.49253$	&	$-0.06220$\\
		  &        & $3$ &  0.14477 &0.31524 &$-0.20649$	&	$-0.02177$\\
	      &        & $4$ &  0.07940 &0.20498 &$-0.12160$	&	$-0.01066$\\
		  &        & $5$ &  0.05095 &0.15421 &$-0.08359$	&	$-0.00614$\\
		  &        & $6$ &  0.03590 &0.12565 &$-0.06286$	&	$-0.00392$\\
				\cline{2-7}
        & $3.5$  & $2$ &  0.43830 &0.67713 &$-0.53995$	&	$-0.05014$\\
		  &        & $3$ &  0.17856 &0.31011 &$-0.23120$	&	$-0.01805$\\
	      &        & $4$ &  0.10157 &0.19801 &$-0.13830$	&	$-0.00906$\\
		  &        & $5$ &  0.06729 &0.14622 &$-0.09616$	&	$-0.00533$\\
		  &        & $6$ &  0.04878 &0.11708 &$-0.07294$	&	$-0.00345$\\
				\cline{2-7}
        & $4.5$  & $2$ &  0.48296 &0.67790 &$-0.56909$	&	$-0.04180$\\
		  &        & $3$ &  0.20117 &0.30848 &$-0.24644$	&	$-0.01529$\\
	      &        & $4$ &  0.11671 &0.19524 &$-0.14865$	&	$-0.00778$\\
		  &        & $5$ &  0.07866 &0.14282 &$-0.10399$	&	$-0.00463$\\
		  &        & $6$ &  0.05788 &0.11331 &$-0.07923$	&	$-0.00303$\\
				\hline
			\end{tabular}
			\end{center}
			\end{table}

\begin{table}[H]
{\bf Table 5:} continued...\vspace{-0.1cm}
\begin{center}
\begin{tabular}{|c|c|c|c|c|c|c|}
				\hline
$k$   & $\theta$ & $n$ &  $V_1$   &  $V_2$  & $V_3$ &$V_4$\\
				\hline
$2$     & $0.75$ & $2$ &  0.22127 &0.69139 &$-0.36184$	&	$-0.09651$\\
		  &        & $3$ &  0.07299 &0.33379 &$-0.13156$	&	$-0.03103$\\
	      &        & $4$ &  0.03304 &0.22973 &$-0.06709$	&	$-0.01373$\\
		  &        & $5$ &  0.01738 &0.18429 &$-0.04041$	&	$-0.00711$\\
		  &        & $6$ &  0.01007 &0.15998 &$-0.02708$	&	$-0.00409$\\
				\cline{2-7}
        & $1.5$  & $2$ &  0.38117 &0.68059 &$-0.49849$	&	$-0.07470$\\
		  &        & $3$ &  0.14642 &0.30610 &$-0.20199$	&	$-0.02705$\\
	      &        & $4$ &  0.07735 &0.19233 &$-0.11334$	&	$-0.01341$\\
		  &        & $5$ &  0.04730 &0.14066 &$-0.07394$	&	$-0.00769$\\
		  &        & $6$ &  0.03162 &0.11230 &$-0.05286$	&	$-0.00484$\\
				\cline{2-7}
        & $2.5$  & $2$ &  0.48740 &0.69160 &$-0.57597$	&	$-0.05453$\\
		  &        & $3$ &  0.19901 &0.30738 &$-0.24309$	&	$-0.02074$\\
	      &        & $4$ &  0.11153 &0.18883 &$-0.14126$	&	$-0.01077$\\
		  &        & $5$ &  0.07210 &0.13419 &$-0.09484$	&	$-0.00644$\\
		  &        & $6$ &  0.05077 &0.10385 &$-0.06940$	&	$-0.00421$\\
				\cline{2-7}
        & $3.5$  & $2$ &  0.54458 &0.69966 &$-0.61473$	&	$-0.04255$\\
		  &        & $3$ &  0.22821 &0.31053 &$-0.26386$	&	$-0.01653$\\
	      &        & $4$ &  0.13109 &0.18970 &$-0.15553$	&	$-0.00876$\\
		  &        & $5$ &  0.08670 &0.13365 &$-0.10565$	&	$-0.00534$\\
		  &        & $6$ &  0.06235 &0.10237 &$-0.07805$	&	$-0.00354$\\
				\cline{2-7}
        & $4.5$  & $2$ &  0.58008 &0.7051  &$-0.63794$	&	$-0.03480$\\
		  &        & $3$ &  0.24661 &0.31302 &$-0.27635$	&	$-0.01369$\\
	      &        & $4$ &  0.14358 &0.19084 &$-0.16415$	&	$-0.00733$\\
		  &        & $5$ &  0.09615 &0.13398 &$-0.11222$	&	$-0.00451$\\
		  &        & $6$ &  0.06993 &0.10214 &$-0.08332$	&	$-0.00302$\\
				\hline
			\end{tabular}
			\end{center}
			\end{table}

\begin{table}[H]
{\bf Table 5:} continued...\vspace{-0.1cm}
\begin{center}
\begin{tabular}{|c|c|c|c|c|c|c|}
				\hline
$k$   & $\theta$ & $n$ &  $V_1$   &  $V_2$  & $V_3$&$V_4$ \\
				\hline

$3$     & $0.75$ & $2$ &  0.30314 &0.67998 &$-0.43325$	&	$-0.09513$\\
		  &        & $3$ &  0.10861 &0.30870 &$-0.16450$	&	$-0.03313$\\
	      &        & $4$ &  0.05284 &0.19790 &$-0.08590$	&	$-0.01564$\\
		  &        & $5$ &  0.02953 &0.14889 &$-0.05210$	&	$-0.00850$\\
		  &        & $6$ &  0.01798 &0.12280 &$-0.03474$	&	$-0.00505$\\
				\cline{2-7}
        & $1.5$  & $2$ &  0.46478 &0.69891 &$-0.56283$	&	$-0.06827$\\
		  &        & $3$ &  0.18665 &0.30846 &$-0.23329$	&	$-0.02620$\\
	      &        & $4$ &  0.10201 &0.18766 &$-0.13218$	&	$-0.01361$\\
		  &        & $5$ &  0.06396 &0.13211 &$-0.08620$	&	$-0.00810$\\
		  &        & $6$ &  0.04353 &0.10148 &$-0.06118$	&	$-0.00523$\\
				\cline{2-7}
        & $2.5$  & $2$ &  0.56164 &0.71809 &$-0.63214$	&	$-0.04804$\\
		  &        & $3$ &  0.23612 &0.31749 &$-0.27102$	&	$-0.01920$\\
	      &        & $4$ &  0.13499 &0.19193 &$-0.15833$	&	$-0.01037$\\
		  &        & $5$ &  0.08836 &0.13336 &$-0.10608$	&	$-0.00640$\\
		  &        & $6$ &  0.06265 &0.10062 &$-0.07706$	&	$-0.00428$\\
				\cline{2-7}
        & $3.5$  & $2$ &  0.61129 &0.72880 &$-0.66588$	&	$-0.03685$\\
		  &        & $3$ &  0.26208 &0.32328 &$-0.28956$	&	$-0.01498$\\
	      &        & $4$ &  0.15268 &0.19549 &$-0.17133$	&	$-0.00823$\\
		  &        & $5$ &  0.10174 &0.13553 &$-0.11606$	&	$-0.00516$\\
		  &        & $6$ &  0.07336 &0.10181 &$-0.08513$	&	$-0.00350$\\
				\cline{2-7}
        & $4.5$  & $2$ &  0.64136 &0.73546 &$-0.68581$	&	$-0.02985$\\
		  &        & $3$ &  0.27797 &0.32704 &$-0.30056$	&	$-0.01225$\\
	      &        & $4$ &  0.16363 &0.19795 &$-0.17907$	&	$-0.00680$\\
		  &        & $5$ &  0.11011 &0.13719 &$-0.12204$	&	$-0.00430$\\
		  &        & $6$ &  0.08012 &0.10291 &$-0.08998$	&	$-0.00294$\\
				\hline
			\end{tabular}
			\end{center}
			\end{table}

\begin{table}[H]
{\bf Table 6:} Simulated quantiles of $T_1$ and $T_2$.\vspace{-0.1cm}
\begin{center}
\resizebox{\textwidth}{!}{
\begin{tabular}{|c|c|c|c|c|c|c|c|c|c|c|}
\hline
 &&&			\multicolumn{4}{c|}{$T_1$} 	&\multicolumn{4}{c|}{$T_2$}\\
\cline{4-11}		
$k$&$\theta$	&$n$	&0.025	&0.05		&0.95		&0.975	&0.025	&0.05		&0.95		&0.975\\
\hline
1&	0.75	&2	&$-0.7673	$&$-0.7180	$&$15.2126	$&$30.6925	$&$-1.0942	$&$-1.0671	$&$2.0643	$&$2.4709$\\
&		&3	&$-0.8982	$&$-0.8328	$&$5.9282	$&$9.4938	$&$-1.3590	$&$-1.2925	$&$1.9171	$&$2.2055$\\
&		&4	&$-0.9758	$&$-0.8908	$&$4.5303	$&$6.8780	$&$-1.4636	$&$-1.3685	$&$1.7941	$&$2.0390$\\
&		&5	&$-1.0449	$&$-0.9345	$&$3.8206	$&$5.3191	$&$-1.4903	$&$-1.3802	$&$1.7655	$&$1.9809$\\
&		&6	&$-1.0579	$&$-0.9479	$&$3.7730	$&$5.2217	$&$-1.5017	$&$-1.3837	$&$1.6857	$&$1.9153$\\
\cline{2-11}
&	1.5	&2	&$-0.8018	$&$-0.7452	$&$16.1365	$&$34.5282	$&$-1.1645	$&$-1.1305	$&$1.9433	$&$2.3739$\\
\hline
			\end{tabular}
	}
			\end{center}
			\end{table}
\begin{table}[H]
{\bf Table 6:} continued...\vspace{-0.1cm}
\begin{center}
\resizebox{\textwidth}{!}{
\begin{tabular}{|c|c|c|c|c|c|c|c|c|c|c|}
\hline
 &&&			\multicolumn{4}{c|}{$T_1$} 	&\multicolumn{4}{c|}{$T_2$}\\
\cline{4-11}	
$k$&$\theta$	&$n$	&0.025	&0.05		&0.95		&0.975	&0.025	&0.05		&0.95		&0.975\\
\hline
2&	0.75	&2	&$-0.7990	$&$-0.7440	$&$16.4358	$&$33.8474	$&$-1.1672	$&$-1.1327	$&$1.9395	$&$2.3697$\\
&		&3	&$-0.9583	$&$-0.8817	$&$6.1137	$&$9.9989	$&$-1.5276	$&$-1.4094	$&$1.8192	$&$2.1227$\\
&		&4	&$-1.0347	$&$-0.9504	$&$4.3106	$&$6.2565	$&$-1.6878	$&$-1.5430	$&$1.6765	$&$1.9223$\\
&		&5	&$-1.1081	$&$-1.0087	$&$3.6158	$&$5.0330	$&$-1.7504	$&$-1.5768	$&$1.6175	$&$1.8473$\\
&		&6	&$-1.1595	$&$-1.0491	$&$3.3388	$&$4.5105	$&$-1.8280	$&$-1.6390	$&$1.5577	$&$1.7701$\\
\cline{2-11}
&	1.5	&2	&$-0.8176	$&$-0.7566	$&$17.2903	$&$36.2768	$&$-1.1747	$&$-1.1364	$&$2.0068	$&$2.4970$\\
\hline
			\end{tabular}
	}
			\end{center}
			\end{table}
\begin{table}[H]
{\bf Table 6:} continued...\vspace{-0.1cm}
\begin{center}
\resizebox{\textwidth}{!}{
\begin{tabular}{|c|c|c|c|c|c|c|c|c|c|c|}
\hline
 &&&			\multicolumn{4}{c|}{$T_1$} 	&\multicolumn{4}{c|}{$T_2$}\\
\cline{4-11}	
$k$&$\theta$	&$n$	&0.025	&0.05		&0.95		&0.975	&0.025	&0.05		&0.95		&0.975\\
\hline
3&	0.75	&2	&$-0.8100	$&$-0.7544	$&$17.1830	$&$35.7986	$&$-1.1750	$&$-1.1386	$&$2.0056	$&$2.4810$\\
&		&3	&$-0.9785	$&$-0.8985	$&$6.2390	$&$9.9324	$&$-1.5546	$&$-1.4270	$&$1.7781	$&$2.1424$\\
&		&4	&$-1.0791	$&$-0.9796	$&$4.0958	$&$6.2641	$&$-1.7360	$&$-1.5727	$&$1.7026	$&$1.9784$\\
&		&5	&$-1.1494	$&$-1.0438	$&$3.6245	$&$5.2769	$&$-1.8766	$&$-1.6442	$&$1.6400	$&$1.8546$\\
&		&6	&$-1.1892	$&$-1.0653	$&$3.2498	$&$4.5584	$&$-1.9333	$&$-1.7088	$&$1.5678	$&$1.7739$\\
\cline{2-11}
&	1.5	&2	&$-0.8118	$&$-0.7463	$&$18.1928	$&$33.9378	$&$-1.1563	$&$-1.1228	$&$2.0054	$&$2.4965$\\
\hline
			\end{tabular}
	}
			\end{center}
			\end{table}

\begin{table}[H]
{\bf Table 7:} Simulated quantiles of $T_3$ and $T_4$.
\begin{center}
\resizebox{\textwidth}{!}{
\begin{tabular}{|c|c|c|c|c|c|c|c|c|c|c|}
\hline
 &&&			\multicolumn{4}{c|}{$T_3$} 	&\multicolumn{4}{c|}{$T_4$}\\
\cline{4-11}	
$k$&$\theta$	&$n$	&0.025	&0.05		&0.95		&0.975	&0.025	&0.05		&0.95		&0.975\\
\hline
1	&0.75	&2	&$-0.8729	$&$-0.7469	$&$39.9467	$&$79.4889	$&$-1.9858	$&$-1.9586	$&$1.1728	$&$1.5794$\\
	&	&3	&$-0.9508	$&$-0.8387	$&$10.7438	$&$16.8522	$&$-2.0265	$&$-1.9599	$&$1.2497	$&$1.5380$\\
	&	&4	&$-1.0262	$&$-0.8979	$&$7.2837	$&$10.8268	$&$-2.0514	$&$-1.9562	$&$1.2063	$&$1.4512$\\
	&	&5	&$-1.1029	$&$-0.9453	$&$5.8435	$&$7.9829	$&$-2.0382	$&$-1.9281	$&$1.2176	$&$1.4331$\\
	&	&6	&$-1.1040	$&$-0.9516	$&$5.5921	$&$7.6001	$&$-2.0251	$&$-1.9070	$&$1.1623	$&$1.3920$\\
\cline{2-11}
	&1.5	&2	&$-0.8835	$&$-0.7328	$&$44.2170	$&$93.1873	$&$-1.9977	$&$-1.9637	$&$1.1102	$&$1.5408$\\
\hline
			\end{tabular}
	}
			\end{center}
			\end{table}
\begin{table}[H]
{\bf Table 7:} continued... 
\begin{center}
\resizebox{\textwidth}{!}{
\begin{tabular}{|c|c|c|c|c|c|c|c|c|c|c|}
\hline
 &&&			\multicolumn{4}{c|}{$T_3$} 	&\multicolumn{4}{c|}{$T_4$}\\
\cline{4-11}	
$k$&$\theta$	&$n$	&0.025	&0.05		&0.95		&0.975	&0.025	&0.05		&0.95		&0.975\\
\hline
2	&0.75	&2	&$-0.8744	$&$-0.7348	$&$42.9065	$&$87.1367	$&$-1.9987	$&$-1.9642	$&$1.1080	$&$1.5382$\\
	&	&3	&$-0.9532	$&$-0.8301	$&$10.4114	$&$16.6549	$&$-2.1053	$&$-1.9872	$&$1.2415	$&$1.5449$\\
	&	&4	&$-1.0105	$&$-0.8945	$&$6.3426	$&$9.0196	$&$-2.1671	$&$-2.0223	$&$1.1972	$&$1.4430$\\
	&	&5	&$-1.0885	$&$-0.9611	$&$4.9665	$&$6.7831	$&$-2.1797	$&$-2.0061	$&$1.1882	$&$1.4180$\\
	&	&6	&$-1.1441	$&$-1.0077	$&$4.4088	$&$5.8550	$&$-2.2279	$&$-2.0390	$&$1.1577	$&$1.3702$\\
\cline{2-11}
	&1.5	&2	&$-0.9176	$&$-0.7516	$&$48.3670	$&$100.0427	$&$-1.9997	$&$-1.9614	$&$1.1818	$&$1.6720$\\
\hline
			\end{tabular}
	}
			\end{center}
			\end{table}
\begin{table}[H]
{\bf Table 7:} continued... 
\begin{center}
\resizebox{\textwidth}{!}{
\begin{tabular}{|c|c|c|c|c|c|c|c|c|c|c|}
\hline
 &&&			\multicolumn{4}{c|}{$T_3$} 	&\multicolumn{4}{c|}{$T_4$}\\
\cline{4-11}	
$k$&$\theta$	&$n$	&0.025	&0.05		&0.95		&0.975	&0.025	&0.05		&0.95		&0.975\\
\hline
3	&0.75	&2	&$-0.8940	$&$-0.7487	$&$46.1976	$&$94.9187	$&$-1.9996	$&$-1.9632	$&$1.1810	$&$1.6564$\\
	&	&3	&$-0.9588	$&$-0.8303	$&$10.6314	$&$16.5626	$&$-2.1102	$&$-1.9826	$&$1.2225	$&$1.5868$\\
	&	&4	&$-1.0364	$&$-0.9021	$&$5.9550	$&$8.8845	$&$-2.1809	$&$-2.0176	$&$1.2578	$&$1.5335$\\
	&	&5	&$-1.1032	$&$-0.9717	$&$4.8443	$&$6.9030	$&$-2.2624	$&$-2.0300	$&$1.2541	$&$1.4687$\\
	&	&6	&$-1.1427	$&$-0.9950	$&$4.1493	$&$5.7094	$&$-2.2837	$&$-2.0592	$&$1.2174	$&$1.4234$\\
\cline{2-11}
	&1.5	&2	&$-0.9193	$&$-0.7345	$&$52.6935	$&$97.1107	$&$-1.9923	$&$-1.9588	$&$1.1694	$&$1.6605$\\
\hline
			\end{tabular}
	}
			\end{center}
			\end{table}

\section{Prediction of Future $k$-Record Values}\label{sec1}
Prediction of future events based on past and present knowledge is a fundamental problem of statistics, arising in many contexts and producing varied solutions. Prediction of future records becomes a problem of great interest. For example, while studying the minimum temperature, having observed the record values until the present time, we will be naturally interested in predicting the minimum temperature that is to be expected when the present record is broken for the first time in the future. Based on the usual records ($k=1$), the problem of predicting future records has been studied by Ahsanullah (1980) and Berred (1998). In this section, we consider prediction for future lower \emph{k}-record values from the UG distribution based on past observed lower \emph{k}-records.
\par

\vspace{0.25cm}
Suppose that the vector of lower \emph{k}-record values $ {\bf R}=(R_{1(k)},R_{2(k)},\ldots,R_{n(k)})$ have been observed from the UG distribution. The problem of interest then is to predict the value of the next lower \emph{k}-record $R_{n+1(k)}$. For a location-scale family with location parameter $\mu$ and scale parameter $\sigma$, the best linear unbiased predictor (BLUP) of $(n+1)$-th  record value was considered by Raqab  (2002). For \emph{k}-record values, the BLUP of $R_{n+1(k)}$ can be written as (see for example Burkschat, 2010)
\begin{equation}\label{4.1}
R^*_{n+1(k)}=\mu^{\ast}+\alpha_{n+1(k)}\sigma ^{\ast}+{\bm \omega}^T \textbf{B}^{-1}({\bf R} -\mu^{\ast}  {\bf 1}-\sigma ^{\ast} \bm{\alpha}),
\end{equation}
where $\mu^{\ast}$ and $\sigma^{\ast}$ are the BLUEs of $\mu$ and $\sigma$   based on the first $n$ lower $k$-record values and  $\alpha_{n+1(k)}$ is the of mean of $(n+1)$-th lower \emph{k}-record value from the standardized UG distribution, $\bm \alpha$ is the vector of  means of the first $n$ lower $k$-records from the standardized UG distribution and
\[
{\bm \omega}^T=\big(Cov(Z_{1(k)},Z_{n+1(k)}), Cov(Z_{2(k)},Z_{n+1(k)}), \cdots, Cov(Z_{n(k)},Z_{n+1(k)})\big),
\]
in which
$Z_{n+1(k)}=\dfrac{R_{n+1(k)}-\mu}{\sigma}$.

\vspace{0.25cm}
Moreover, the mean squared prediction error (MSPE) of
$R^*_{n+1(k)}$ is given by (Burkschat, 2010)
\begin{eqnarray*}
{\rm MSPE}(R^*_{n+1(k)})&=&E[(R^*_{n+1(k)}-R_{n+1(k)})^2]\\&=&
\sigma^2[(1-{\bm  \omega}^T\textbf{B}^{-1}\mathbf{1})^2V_1+(\alpha_{n+1(k)}-{\bm  \omega}^T\textbf{B}^{-1}{\bm  \alpha})^2V_2-{\bm \omega}^T\textbf{B}^{-1}{\bm  \omega}\\&&+2(1-{\bm  \omega}^T\textbf{B}^{-1}\mathbf{1})(\alpha_{n+1(k)}-{\bm  \omega}^T\textbf{B}^{-1}{\bm  \alpha})V_3+Var(Z_{n+1(k)})],
\end{eqnarray*}
where
$\alpha_{n+1(k)}=E(Z_{n+1(k)})$.

\vspace{0.25cm}
Let us now consider the best linear invariant predictor (BLIP) of
the next upper $k$-record value. From the results of Mann (1969), the
BLIP of $R_{n+1(k)}$ can be obtained based on the BLUP of $R_{n+1(k)}$ as follows (see
also Arnold et al., 1998, p. 153)
\[\widetilde{R}_{n+1(k)}=R^*_{n+1(k)}-\left(\frac{V_4}{1+V_2}\right) \sigma^*,\]
where $R^*_{n+1(k)}$ is the BLUP of $R_{n+1(k)}$ and
\[
V_4=(1-{\bm  \omega}^T\textbf{B}^{-1}\mathbf{1})V_3+(\alpha_{n+1(k)}-{\bm  \omega}^T\textbf{B}^{-1}{\bm  \alpha})V_2.
\]
In Table 5, we reported the values of $V_4$ for different values of $n$ and
$\theta$.

\vspace{0.25cm}
The MSPE of $\widetilde{R}_{n+1(k)}$ is given by (Burkschat, 2010)
\begin{eqnarray*}
{\rm MSPE}(\widetilde{R}_{n+1(k)})&=&E[(\widetilde{R}_{n+1(k)}-R_{n+1(k)})^2]\\&=&
\sigma^2\bigg{[}\frac{{\bm  \alpha}^T\textbf{B}^{-1}{\bm \alpha}+1}{\Delta}(1-{\bm  \omega}^T\textbf{B}^{-1}\mathbf{1})^2+\frac{\textbf{1}^T\textbf{B}^{-1}\textbf{1}}{\Delta}
(\alpha_{n+1(k)}-{\bm  \omega}^T\textbf{B}^{-1}{\bm  \alpha})^2\\&&\quad \quad -{\bm  \omega}^T\textbf{B}^{-1}{\bm  \omega}  -2
\frac{{\bm \alpha}^T\textbf{B}^{-1}\textbf{1}}{\Delta}(1-{\bm  \omega}^T\textbf{B}^{-1}\mathbf{1})(\alpha_{n+1(k)}-{\bm  \omega}^T\textbf{B}^{-1}{\bm  \alpha})
\\&& \quad\quad +Var(Z_{n+1(k)})\bigg{]},
\end{eqnarray*}
where
\begin{eqnarray*}
\Delta=({\bm \alpha}^T\textbf{B}^{-1}{\bm  \alpha}+1)(\textbf{1}^T\textbf{B}^{-1}\textbf{1})-({\bm \alpha}^T\textbf{B}^{-1}\textbf{1})^2.
\end{eqnarray*}

Now we compare the BLUP and BLIP of $R_{n+1(k)}$ using the REC. The REC of
$\widetilde{R}_{n+1(k)}$ relative to $R^*_{n+1(k)}$ is
\[REC(\widetilde{R}_{n+1(k)},R^*_{n+1(k)})=\frac{{\rm MSPE}(R^*_{n+1(k)})}{{\rm MSPE}(\widetilde{R}_{n+1(k)})}.\]
In Table 8, we presented the REC values of $\widetilde{R}_{n+1(k)}$ relative to
$R^*_{n+1(k)}$ for different choices of $n$ and $\theta$. From
Table 8, we observe that the BLIP works better than BLUP in terms of
MSPE.

\vspace{0.25cm}
Suppose we are now interested in prediction intervals (PIs) for
$R_{n+1(k)}$. The PIs  can be constructed  using the pivotal
quantities (Balakrishnan and Chan, 1998)
\[T^* _{1}=\frac{R_{n+1(k)}-R_{n(k)}}{\sigma^*},\]
and
\[T^* _{2}=\frac{R_{n+1(k)}-R_{n(k)}}{\widetilde{\sigma}}.\]
Constructing such PIs requires the percentage points of $T^*_1$ and
$T^*_2$. In Table 9, we presented the simulated percentage points of
$T^*_1$ and $T^*_2$ using a Monte Carlo method based on 10000 runs and
different choices of $n$ and  $\theta$. Using the pivotal quantity
$T^*_1$, a $100(1-\gamma)\%$ PI for $R_{n+1(k)}$ is given by
\begin{equation}\label{cimu3}
P \bigg{(}R_{n(k)}+\sigma^* T^*_1(\gamma/2)\leq R_{n+1(k)}
\leq R_{n(k)}+\sigma^* T^*_1(1-\gamma/2)\bigg{)}=1-\gamma,
\end{equation}
where $T_1^*(\tau)$ is the left percentage point of $T_1^*$ at $\tau$, i.e. $P(T^*_1<T^*_1(\tau))=\tau$.

\vspace{0.5cm}
Similarly, using the pivotal quantity $T^*_2$, a $100(1-\gamma)\%$
PI for $R_{n+1(k)}$ is given by
\begin{equation}\label{cimu4}
P \bigg{(}R_{n(k)}+\widetilde{\sigma} T^*_2(\gamma/2)\leq R_{n+1(k)}
\leq R_{n(k)}+\widetilde{\sigma} T^*_2(1-\gamma/2)\bigg{)}=1-\gamma,
\end{equation}
where $T_2^*(\tau)$ is the left percentage point of $T_2^*$ at $\tau$, i.e. $P(T^*_2<T^*_2(\tau))=\tau$.

\begin{table}[H]
{\bf Table 8:} The REC values of $\widetilde{R}_{n+1(k)}$ relative to $R^*_{n+1(k)}$.\vspace{-0.1cm}
\begin{center}
\begin{tabular}{|c|c|c|c|c|c|c|}
\hline
$k$&	$n$&	$\theta=0.75$&	$\theta=1.5$&	$\theta=2.5$&	$\theta=3.5$&	$\theta=4.5$\\
\hline
1&	2&	1.24169&	1.26678&	1.27618&	1.28014&	1.28233\\
&	3&	1.07436&	1.08672&	1.09151&	1.09358&	1.09474\\
&	4&	1.03482&	1.04132&	1.04398&	1.04518&	1.04587\\
&	5&	1.02065&	1.02390&	1.02538&	1.02608&	1.02650\\
&	6&	1.01432&	1.01566&	1.01642&	1.01683&	1.01709\\
\hline
2&	2&	1.24169&	1.26678&	1.27618&	1.28014&	1.28233\\
&	3&	1.07436&	1.08672&	1.09151&	1.09358&	1.09474\\
&	4&	1.03482&	1.04132&	1.04398&	1.04518&	1.04587\\
&	5&	1.02065&	1.02390&	1.02538&	1.02608&	1.02650\\
&	6&	1.01432&	1.01566&	1.01642&	1.01683&	1.01709\\
\hline
3&	2&	1.24169&	1.26678&	1.27618&	1.28014&	1.28233\\
&	3&	1.07436&	1.08672&	1.09151&	1.09358&	1.09474\\
&	4&	1.03482&	1.04132&	1.04398&	1.04518&	1.04587\\
&	5&	1.02065&	1.02390&	1.02538&	1.02608&	1.02650\\
&	6&	1.01432&	1.01566&	1.01642&	1.01683&	1.01709\\

\hline
			\end{tabular}
			\end{center}
			\end{table}

\begin{table}[H]
{\bf Table 9:} Simulated quantiles of $T_1^\ast$ and $T_2^\ast$.\vspace{-0.5cm}
\begin{center}
\resizebox{\textwidth}{!}{
\begin{tabular}{|c|c|c|c|c|c|c|c|c|c|c|}
\hline
 &&&			\multicolumn{4}{c|}{$T_1^\ast$} 	&\multicolumn{4}{c|}{$T_2^\ast$}\\
\cline{4-11}			
$k$&$\theta$	&$n$	&0.025	&0.05		&0.95		&0.975	&0.025	&0.05		&0.95		&0.975\\
\hline
1&	0.75&	2	&$-4.2883	$&$-2.0201	$&$-0.0040	$&$-0.0019	$&$-7.6969	$&$-3.6259	$&$-0.0071	$&$-0.0033$\\
&&		3	&$-0.6562	$&$-0.4205	$&$-0.0024	$&$-0.0011	$&$-0.9485	$&$-0.6078	$&$-0.0034	$&$-0.0016$\\
&&		4	&$-0.2932	$&$-0.1951	$&$-0.0016	$&$-0.0008	$&$-0.3945	$&$-0.2625	$&$-0.0022	$&$-0.0010$\\
&&		5	&$-0.1607	$&$-0.1166	$&$-0.0012	$&$-0.0006	$&$-0.2089	$&$-0.1516	$&$-0.0016	$&$-0.0008$\\
&&		6	&$-0.1085	$&$-0.0795	$&$-0.0009	$&$-0.0004	$&$-0.1382	$&$-0.1013	$&$-0.0011	$&$-0.0005$\\
\cline{2-11}
&	1.5&	2	&$-3.4072	$&$-1.8017	$&$-0.0046	$&$-0.0022	$&$-5.7723	$&$-3.0524	$&$-0.0079	$&$-0.0038$\\
&&		3	&$-0.6550	$&$-0.4225	$&$-0.0031	$&$-0.0016	$&$-0.8759	$&$-0.5651	$&$-0.0042	$&$-0.0021$\\
&&		4	&$-0.3083	$&$-0.2252	$&$-0.0024	$&$-0.0012	$&$-0.3795	$&$-0.2772	$&$-0.0030	$&$-0.0015$\\
&&		5	&$-0.2152	$&$-0.1580	$&$-0.0020	$&$-0.0010	$&$-0.2544	$&$-0.1867	$&$-0.0023	$&$-0.0012$\\
&&		6	&$-0.1605	$&$-0.1198	$&$-0.0015	$&$-0.0008	$&$-0.1852	$&$-0.1383	$&$-0.0017	$&$-0.0009$\\
\cline{2-11}
&	2.5&	2	&$-2.8387	$&$-1.3615	$&$-0.0037	$&$-0.0020	$&$-4.7647	$&$-2.2853	$&$-0.0063	$&$-0.0033$\\
&&		3	&$-0.5648	$&$-0.3643	$&$-0.0028	$&$-0.0014	$&$-0.7428	$&$-0.4791	$&$-0.0036	$&$-0.0018$\\
&&		4	&$-0.2746	$&$-0.1991	$&$-0.0021	$&$-0.0010	$&$-0.3309	$&$-0.2399	$&$-0.0025	$&$-0.0012$\\
&&		5	&$-0.2012	$&$-0.1487	$&$-0.0016	$&$-0.0008	$&$-0.2323	$&$-0.1716	$&$-0.0019	$&$-0.0010$\\
&&		6	&$-0.1463	$&$-0.1111	$&$-0.0015	$&$-0.0008	$&$-0.1646	$&$-0.1250	$&$-0.0017	$&$-0.0009$\\
\cline{2-11}
&	3.5&	2	&$-2.2998	$&$-1.1448	$&$-0.0032	$&$-0.0016	$&$-3.8571	$&$-1.9200	$&$-0.0053	$&$-0.0027$\\
&&		3	&$-0.4710	$&$-0.3012	$&$-0.0023	$&$-0.0011	$&$-0.6170	$&$-0.3945	$&$-0.0030	$&$-0.0014$\\
&&		4	&$-0.2428	$&$-0.1725	$&$-0.0019	$&$-0.0009	$&$-0.2908	$&$-0.2066	$&$-0.0023	$&$-0.0011$\\
&&		5	&$-0.1722	$&$-0.1253	$&$-0.0015	$&$-0.0007	$&$-0.1974	$&$-0.1436	$&$-0.0017	$&$-0.0008$\\
&&		6	&$-0.1310	$&$-0.0978	$&$-0.0013	$&$-0.0007	$&$-0.1463	$&$-0.1092	$&$-0.0015	$&$-0.0008$\\
\cline{2-11}
&	4.5&	2	&$-1.7988	$&$-0.9013	$&$-0.0025	$&$-0.0013	$&$-3.0183	$&$-1.5124	$&$-0.0043	$&$-0.0021$\\
&&		3	&$-0.4050	$&$-0.2605	$&$-0.0018	$&$-0.0009	$&$-0.5299	$&$-0.3409	$&$-0.0024	$&$-0.0011$\\
&&		4	&$-0.2026	$&$-0.1463	$&$-0.0016	$&$-0.0008	$&$-0.2422	$&$-0.1749	$&$-0.0019	$&$-0.0010$\\
&&		5	&$-0.1494	$&$-0.1085	$&$-0.0013	$&$-0.0006	$&$-0.1708	$&$-0.1240	$&$-0.0015	$&$-0.0007$\\
&&		6	&$-0.1161	$&$-0.0874	$&$-0.0012	$&$-0.0006	$&$-0.1293	$&$-0.0973	$&$-0.0014	$&$-0.0007$\\
\hline
			\end{tabular}
	}
			\end{center}
			\end{table}
\begin{table}[H]
{\bf Table 9:} continued...\vspace{-0.5cm}
\begin{center}
\resizebox{\textwidth}{!}{
\begin{tabular}{|c|c|c|c|c|c|c|c|c|c|c|}
\hline
 &&&			\multicolumn{4}{c|}{$T_1^\ast$} 	&\multicolumn{4}{c|}{$T_2^\ast$}\\
\cline{4-11}	
$k$&$\theta$	&$n$		&0.025	&0.05		&0.95		&0.975	&0.025	&0.05		&0.95		&0.975\\
\hline
2&	0.75&	2	&$-4.3459	$&$-2.1838	$&$-0.0050	$&$-0.0025	$&$-7.3507	$&$-3.6937	$&$-0.0084	$&$-0.0042$\\
&&		3	&$-0.7805	$&$-0.5290	$&$-0.0031	$&$-0.0015	$&$-1.0410	$&$-0.7056	$&$-0.0042	$&$-0.0020$\\
&&		4	&$-0.3763	$&$-0.2634	$&$-0.0023	$&$-0.0011	$&$-0.4627	$&$-0.3239	$&$-0.0028	$&$-0.0013$\\
&&		5	&$-0.2337	$&$-0.1676	$&$-0.0017	$&$-0.0008	$&$-0.2768	$&$-0.1985	$&$-0.0020	$&$-0.0010$\\
&&		6	&$-0.1585	$&$-0.1163	$&$-0.0012	$&$-0.0006	$&$-0.1839	$&$-0.1349	$&$-0.0014	$&$-0.0007$\\
\cline{2-11}
&	1.5&	2	&$-3.5896	$&$-1.6659	$&$-0.0045	$&$-0.0023	$&$-6.0326	$&$-2.7996	$&$-0.0076	$&$-0.0039$\\
&&		3	&$-0.6905	$&$-0.4389	$&$-0.0032	$&$-0.0015	$&$-0.9018	$&$-0.5733	$&$-0.0042	$&$-0.0019$\\
&&		4	&$-0.3522	$&$-0.2511	$&$-0.0025	$&$-0.0011	$&$-0.4199	$&$-0.2993	$&$-0.0029	$&$-0.0013$\\
&&		5	&$-0.2359	$&$-0.1714	$&$-0.0019	$&$-0.0009	$&$-0.2690	$&$-0.1955	$&$-0.0022	$&$-0.0010$\\
&&		6	&$-0.1712	$&$-0.1308	$&$-0.0015	$&$-0.0007	$&$-0.1904	$&$-0.1455	$&$-0.0017	$&$-0.0007$\\
\cline{2-11}
&	2.5&	2	&$-2.5916	$&$-1.2622	$&$-0.0030	$&$-0.0016	$&$-4.3840	$&$-2.1352	$&$-0.0051	$&$-0.0026$\\
&&		3	&$-0.5198	$&$-0.3496	$&$-0.0026	$&$-0.0013	$&$-0.6796	$&$-0.4570	$&$-0.0034	$&$-0.0017$\\
&&		4	&$-0.2778	$&$-0.2042	$&$-0.0019	$&$-0.0009	$&$-0.3303	$&$-0.2427	$&$-0.0023	$&$-0.0011$\\
&&		5	&$-0.1899	$&$-0.1466	$&$-0.0016	$&$-0.0008	$&$-0.2154	$&$-0.1662	$&$-0.0019	$&$-0.0009$\\
&&		6	&$-0.1501	$&$-0.1133	$&$-0.0014	$&$-0.0007	$&$-0.1657	$&$-0.1251	$&$-0.0016	$&$-0.0008$\\
\cline{2-11}
&	3.5&	2	&$-1.8265	$&$-0.8932	$&$-0.0026	$&$-0.0013	$&$-3.1045	$&$-1.5181	$&$-0.0045	$&$-0.0022$\\
&&		3	&$-0.3890	$&$-0.2581	$&$-0.0021	$&$-0.0010	$&$-0.5098	$&$-0.3382	$&$-0.0027	$&$-0.0013$\\
&&		4	&$-0.2275	$&$-0.1630	$&$-0.0017	$&$-0.0008	$&$-0.2706	$&$-0.1939	$&$-0.0020	$&$-0.0009$\\
&&		5	&$-0.1563	$&$-0.1162	$&$-0.0012	$&$-0.0007	$&$-0.1772	$&$-0.1318	$&$-0.0014	$&$-0.0008$\\
&&		6	&$-0.1248	$&$-0.0934	$&$-0.0012	$&$-0.0005	$&$-0.1375	$&$-0.1030	$&$-0.0013	$&$-0.0006$\\
\cline{2-11}
&	4.5&	2	&$-1.4966	$&$-0.7438	$&$-0.0022	$&$-0.0011	$&$-2.5518	$&$-1.2682	$&$-0.0038	$&$-0.0019$\\
&&		3	&$-0.3667	$&$-0.2244	$&$-0.0017	$&$-0.0007	$&$-0.4815	$&$-0.2947	$&$-0.0022	$&$-0.0010$\\
&&		4	&$-0.1971	$&$-0.1420	$&$-0.0014	$&$-0.0007	$&$-0.2348	$&$-0.1691	$&$-0.0017	$&$-0.0008$\\
&&		5	&$-0.1346	$&$-0.0991	$&$-0.0012	$&$-0.0006	$&$-0.1527	$&$-0.1123	$&$-0.0014	$&$-0.0007$\\
&&		6	&$-0.1077	$&$-0.0786	$&$-0.0010	$&$-0.0005	$&$-0.1187	$&$-0.0866	$&$-0.0011	$&$-0.0005$\\
\hline
			\end{tabular}
	}
			\end{center}
			\end{table}
\begin{table}[H]
{\bf Table 9:} continued...\vspace{-0.5cm}
\begin{center}
\resizebox{\textwidth}{!}{
\begin{tabular}{|c|c|c|c|c|c|c|c|c|c|c|}
\hline
 &&&			\multicolumn{4}{c|}{$T_1^\ast$} 	&\multicolumn{4}{c|}{$T_2^\ast$}\\
\cline{4-11}	
$k$&$\theta$	&$n$	&0.025	&0.05		&0.95		&0.975	&0.025	&0.05		&0.95		&0.975\\

\hline
3&	0.75&	2	&$-3.8129	$&$-2.0512	$&$-0.0058	$&$-0.0028	$&$-6.4057	$&$-3.4460	$&$-0.0097	$&$-0.0047$\\
&&		3	&$-0.8329	$&$-0.5563	$&$-0.0042	$&$-0.0022	$&$-1.0900	$&$-0.7280	$&$-0.0055	$&$-0.0029$\\
&&		4	&$-0.4191	$&$-0.3075	$&$-0.0026	$&$-0.0013	$&$-0.5020	$&$-0.3684	$&$-0.0032	$&$-0.0015$\\
&&		5	&$-0.2880	$&$-0.2009	$&$-0.0019	$&$-0.0010	$&$-0.3308	$&$-0.2308	$&$-0.0022	$&$-0.0011$\\
&&		6	&$-0.1819	$&$-0.1338	$&$-0.0016	$&$-0.0008	$&$-0.2042	$&$-0.1503	$&$-0.0018	$&$-0.0009$\\
\cline{2-11}
&	1.5&	2	&$-3.2888	$&$-1.5978	$&$-0.0041	$&$-0.0019	$&$-5.5874	$&$-2.7146	$&$-0.0070	$&$-0.0033$\\
&&		3	&$-0.6753	$&$-0.4108	$&$-0.0032	$&$-0.0018	$&$-0.8836	$&$-0.5375	$&$-0.0042	$&$-0.0023$\\
&&		4	&$-0.3665	$&$-0.2564	$&$-0.0025	$&$-0.0012	$&$-0.4353	$&$-0.3045	$&$-0.0030	$&$-0.0015$\\
&&		5	&$-0.2572	$&$-0.1906	$&$-0.0020	$&$-0.0010	$&$-0.2912	$&$-0.2158	$&$-0.0023	$&$-0.0011$\\
&&		6	&$-0.1954	$&$-0.1400	$&$-0.0017	$&$-0.0008	$&$-0.2152	$&$-0.1542	$&$-0.0018	$&$-0.0009$\\
\cline{2-11}
&	2.5&	2	&$-2.0008	$&$-0.9970	$&$-0.0029	$&$-0.0013	$&$-3.4376	$&$-1.7129	$&$-0.0050	$&$-0.0022$\\
&&		3	&$-0.4617	$&$-0.3100	$&$-0.0023	$&$-0.0012	$&$-0.6082	$&$-0.4084	$&$-0.0031	$&$-0.0015$\\
&&		4	&$-0.2642	$&$-0.1825	$&$-0.0018	$&$-0.0009	$&$-0.3149	$&$-0.2175	$&$-0.0022	$&$-0.0011$\\
&&		5	&$-0.1857	$&$-0.1375	$&$-0.0016	$&$-0.0008	$&$-0.2105	$&$-0.1558	$&$-0.0018	$&$-0.0009$\\
&&		6	&$-0.1487	$&$-0.1137	$&$-0.0014	$&$-0.0007	$&$-0.1637	$&$-0.1251	$&$-0.0015	$&$-0.0007$\\
\cline{2-11}
&	3.5&	2	&$-1.7763	$&$-0.8347	$&$-0.0023	$&$-0.0011	$&$-3.0709	$&$-1.4430	$&$-0.0040	$&$-0.0020$\\
&&		3	&$-0.3819	$&$-0.2384	$&$-0.0019	$&$-0.0008	$&$-0.5053	$&$-0.3154	$&$-0.0025	$&$-0.0011$\\
&&		4	&$-0.2148	$&$-0.1475	$&$-0.0015	$&$-0.0007	$&$-0.2568	$&$-0.1763	$&$-0.0018	$&$-0.0009$\\
&&		5	&$-0.1569	$&$-0.1128	$&$-0.0013	$&$-0.0007	$&$-0.1782	$&$-0.1281	$&$-0.0015	$&$-0.0008$\\
&&		6	&$-0.1182	$&$-0.0903	$&$-0.0011	$&$-0.0006	$&$-0.1303	$&$-0.0995	$&$-0.0012	$&$-0.0006$\\
\cline{2-11}
&	4.5&	2	&$-1.3292	$&$-0.6353	$&$-0.0019	$&$-0.0009	$&$-2.3068	$&$-1.1026	$&$-0.0032	$&$-0.0016$\\
&&		3	&$-0.2760	$&$-0.1862	$&$-0.0015	$&$-0.0007	$&$-0.3662	$&$-0.2471	$&$-0.0019	$&$-0.0010$\\
&&		4	&$-0.1760	$&$-0.1261	$&$-0.0012	$&$-0.0005	$&$-0.2109	$&$-0.1510	$&$-0.0014	$&$-0.0006$\\
&&		5	&$-0.1297	$&$-0.0912	$&$-0.0010	$&$-0.0005	$&$-0.1474	$&$-0.1038	$&$-0.0011	$&$-0.0006$\\
&&		6	&$-0.1001	$&$-0.0742	$&$-0.0009	$&$-0.0005	$&$-0.1104	$&$-0.0818	$&$-0.0010	$&$-0.0005$\\
\hline
			\end{tabular}
	}
			\end{center}
			\end{table}
\section{A Simulation Study}
In this section, we conducted a simulation study to evaluate the estimators and predictors proposed in this paper. To this end, we generated randomly $N=2000$ $k$-records from the standardized unit-Gompertz distribution for $k=1,2,3$ and  different choices of $\theta$ and $n$.
In each iteration, the BLUEs and BLIEs of $\mu$ and $\sigma$ as well as the BLUP and BLIP of the $(n+1)$-th $k$-record were calculted. Besides, using the simulated quantiles given in Tables 6, 7 and 9 and the relations  (\ref{cimu}),  (\ref{cisigma}),  (\ref{cimu2}),  (\ref{cisigma2}),  (\ref{cimu3}) and  (\ref{cimu4}), we computed
the $95\%$ CIs for $\mu$ based on the pivotal quantities $T_1$ and $T_3$,  the $95\%$ CIs for $\sigma$ based on the pivotal quantities $T_2$ and $T_4$
 and $95\%$ PIs for $R_{n+1(k)}$ based on  the pivotal quantities $T_1^*$ and $T_2^*$.

\vspace{0.25cm}
We compare the performance of the point estimators and point predictors  using the estimated bias (EB)  and the estimated mean squared error (EMSE).
Let $\hat{\mu}$ be an estimator of $\mu$ and $\hat{\mu}_i$ be the corresponding estimate obtained in the $i$-th iteration of the simulation. Then the EB and the EMSE of $\hat{\mu}$ are given by
\begin{equation}\label{eb}
{\rm EB}(\hat{\mu})=\frac{1}{N}\sum_{i=1}^N (\hat{\mu}_i-\mu),
\end{equation}
and
\begin{equation}\label{emse}
{\rm EMSE}(\hat{\mu})=\frac{1}{N}\sum_{i=1}^N (\hat{\mu}_i-\mu)^2,
\end{equation}
respectively.

\vspace{0.25cm}
Similarly, we can define the EB and EMSE of an estimator of $\sigma$. We calculated the EBs and EMSEs of the BLUEs and BLIEs of $\mu$ and $\sigma$ and the results are given in Table 10.

\vspace{0.25cm}
Besides, we can define the EB and the estimated mean squared prediction error (EMSPE) of a predictor of $R_{n+1(k)}$ similarly to (\ref{eb}) and (\ref{emse}), respectively.
So we can evaluate the point predictors of $R_{n+1(k)}$ using the EB and EMSPE criteria.
We also compare the performance of the $95\%$ CIs and PIs using the average length (AL) and coverage probability (CP).
The results for $95\%$ CIs of the location and scale parameters are given in Table 11 and the simulation results related to the point and interval prediction are reported in Table 12.

\begin{table}[H]
{\bf Table 10:} The EBs and EMSEs of the BLUEs and BLIEs of $\mu$ and $\sigma$.\vspace{-0.1cm}
\begin{center}
\resizebox{\textwidth}{!}{
\begin{tabular}{|c|c|c|c|c|c|c|c|c|c|c|}
\hline
 &&&			\multicolumn{2}{c|}{$\mu^\ast$} 	&\multicolumn{2}{c|}{$\sigma^\ast$}
&			\multicolumn{2}{c|}{$\widetilde{\mu}$} 	&\multicolumn{2}{c|}{$\widetilde{\sigma}$}\\
\cline{4-11}	
$k$&$\theta$	&$n$	&EB	&EMSE		&EB		&EMSE	&EB	&EMSE		&EB		&EMSE\\

\hline
1&	$0.75$&	$2$	&$0.0025$ &	$0.1130$ & 	$-0.0106$& 	$0.8069$&	$0.1460$ &	$0.0760$&	$-0.4488$&	$0.4518$
\\
&&		$3$	&$0.0052$&	$0.0320$ & 	$-0.0249$& 	$0.4317$& 	$0.0654$&	$0.0270$&	$-0.3255$&	$0.3122$
 \\
&&		$4$	&$-0.0030$& 	$0.0128$& 	$0.0350$& 	$0.3575$& $0.0316$&		$0.0110$&	$-0.2307$&	$0.2501$
 \\
&&		$5$	&$0.0024$& 	$0.0067$& $0.0435$& 	$0.2984$& 	$0.0245$&	$0.0061$&	$-0.1974$&	$0.2144$
 \\
&&		$6$	&$0.0041$& 	$0.0036$& 	$0.0493$& 	$0.2758$& 	$0.0197$&	$0.0034$&	$-0.1763$&	$0.1995$
 \\
\cline{2-11}
&	$1.5$&	$2$	&$0.0142$	 & $0.2551$& 	$-0.0189$& 	$0.6834$&$0.2473$ &		0.1654&	$-0.4209$&	0.4152
\\
\hline
			\end{tabular}
	}
			\end{center}
			\end{table}

\begin{table}[H]
{\bf Table 10:} continued...\vspace{-0.1cm}
\begin{center}
\resizebox{\textwidth}{!}{
\begin{tabular}{|c|c|c|c|c|c|c|c|c|c|c|}
\hline
 &&&			\multicolumn{2}{c|}{$\mu^\ast$} 	&\multicolumn{2}{c|}{$\sigma^\ast$}
&			\multicolumn{2}{c|}{$\widetilde{\mu}$} 	&\multicolumn{2}{c|}{$\widetilde{\sigma}$}\\
\cline{4-11}	
$k$&$\theta$	&$n$	&EB	&EMSE		&EB		&EMSE	&EB	&EMSE		&EB		&EMSE\\

\hline
2&	0.75&	2	&0.0272 &	0.2099 & 	$-0.0572$&	0.6612 &	0.2289 &	0.1457&	$-0.4426$&	0.4259
\\
&&		3	&0.0065 & 	0.0722& $-0.0158$&		0.3353 & 	0.1036&	0.0606&	$-0.2621$&	0.2570
\\
&&		4	&$-0.0021$ & 	0.0316& 	0.0053& 	0.2230& 0.0528&		0.0281&	$-0.1825$&	0.1807
\\
&&		5	& $-0.0022$& 	0.0175& 	0.0046&	0.1861 & 	0.0321&	0.0159&	$-0.1517$&	0.1557
\\
&&		6	& $-0.0045$& 	0.0100& 0.0071&		0.1568 &	0.0190 &	0.0091&	$-0.1318$&	0.1338
\\
\cline{2-11}
&	1.5&	2	&$-0.0283$ & 	0.4153& 	0.0456& 	0.7380& 	0.2818&	0.2375&	$-0.3778$&	0.4033\\
\hline
			\end{tabular}
	}
			\end{center}
			\end{table}

\begin{table}[H]
{\bf Table 10:} continued...\vspace{-0.1cm}
\begin{center}
\resizebox{\textwidth}{!}{
\begin{tabular}{|c|c|c|c|c|c|c|c|c|c|c|}
\hline
 &&&			\multicolumn{2}{c|}{$\mu^\ast$} 	&\multicolumn{2}{c|}{$\sigma^\ast$}
&			\multicolumn{2}{c|}{$\widetilde{\mu}$} 	&\multicolumn{2}{c|}{$\widetilde{\sigma}$}\\
\cline{4-11}	
$k$&$\theta$	&$n$	&EB	&EMSE		&EB		&EMSE	&EB	&EMSE		&EB		&EMSE\\

\hline
3&	0.75&	2	& 0.0016& 	0.2902& 	$-0.0065$	& 0.6570& 	0.2578&	0.1858&	$-0.4086$&	0.3997\\
&&		3	&$-0.0030$& 	0.1066& 	0.0061&	0.3062 &	0.1235 &	0.0860&	$-0.2312$&	0.2322 \\
&&		4	&$-0.0056$& 	0.0508& 	0.0198& 	0.1904&	0.0675 &	0.0444&	$-0.1487$&	0.1545\\
&&		5	& $-0.0053$& 	0.0292& $-0.0020$& 	0.1584& 0.0400&	0.0262&	$-0.1313$&	0.1373 \\
&&		6	&0.0008	& 0.0174& 	$-0.0068$	& 0.1259& 	0.0315&	0.0163&	$-0.1155$&	0.1132 \\
\cline{2-11}
&	1.5&	2	&0.0075 & 	0.4662& 	$-0.0110$& 	0.7017& 0.3351&	0.2813&	$-0.4178$&	0.4177\\
\hline
			\end{tabular}
	}
			\end{center}
			\end{table}

\begin{table}[H]
{\bf Table 11:} The ALs and CPs of $95\%$ CIs for $\mu$ and $\sigma$ based on $T_1, T_2, T_3$ and $T_4$. 
\begin{center}
\resizebox{\textwidth}{!}{
\begin{tabular}{|c|c|c|c|c|c|c|c|c|c|c|}
\hline
 &&&			\multicolumn{4}{c|}{CIs for $\mu$} 	&\multicolumn{4}{c|}{CIs for $\sigma$} \\
 \cline{4-11}
 &&&			\multicolumn{2}{c|}{$T_1$} 	&\multicolumn{2}{c|}{$T_3$}
&			\multicolumn{2}{c|}{$T_2$} 	&\multicolumn{2}{c|}{$T_4$}\\
\cline{4-11}	
$k$&$\theta$	&$n$	&AL	&CP		&AL		&CP	&AL	&CP		&AL		&CP\\

\hline
1&	0.75&	2	&10.6032&1.0000&10.6032 &1.0000 &40.2064 &1.0000&40.2058&1.0000\\
&&		3	&1.8229&1.0000&1.8229&1.0000 &10.1022&1.0000&10.10225&1.0000\\
&&		4	& 0.9184&1.0000&0.9184&1.0000 &6.9354&1.0000&6.9354&1.0000\\
&&		5	& 0.5160&1.0000&0.5160&1.0000 &5.1871&1.0000&5.1871&1.0000\\
&&		6	& 0.3745&1.0000&0.3745&1.0000 &4.3775&1.0000&4.377&1.0000\\
\hline
			\end{tabular}
	}
			\end{center}
			\end{table}

\begin{table}[H]
{\bf Table 11:} continued...\vspace{-0.5cm}
\begin{center}
\resizebox{\textwidth}{!}{
\begin{tabular}{|c|c|c|c|c|c|c|c|c|c|c|}
\hline
 &&&			\multicolumn{4}{c|}{CIs for $\mu$} 	&\multicolumn{4}{c|}{CIs for $\sigma$} \\
 \cline{4-11}
 &&&			\multicolumn{2}{c|}{$T_1$} 	&\multicolumn{2}{c|}{$T_3$}
&			\multicolumn{2}{c|}{$T_2$} 	&\multicolumn{2}{c|}{$T_4$}\\
\cline{4-11}	
$k$&$\theta$	&$n$	&AL	&CP		&AL		&CP	&AL	&CP		&AL		&CP\\

\hline
2&	0.75&	2	&15.3656 &1.0000 &15.3656 & 1.0000&31.6666 &1.0000&31.6660&1.0000\\
&&		3	&2.9134 &1.0000 &2.9134 &1.0000 &7.9371 &1.0000&7.9371&1.0000\\
&&		4	&1.3322 & 1.0000&1.3322 &1.0000 &4.7396 &1.0000&4.7396&1.0000\\
&&		5	&0.8133 &1.0000 &0.8133 &1.0000 &3.4814 &1.0000&3.4814&1.0000\\
&&		6	&0.5729 &1.0000 &0.5729 & 1.0000&3.1561 &1.0000&3.1561&1.0000\\
\hline
			\end{tabular}
	}
			\end{center}
			\end{table}

\begin{table}[H]
{\bf Table 11:} continued...\vspace{-0.5cm}
\begin{center}
\resizebox{\textwidth}{!}{
\begin{tabular}{|c|c|c|c|c|c|c|c|c|c|c|}
\hline
 &&&			\multicolumn{4}{c|}{CIs for $\mu$} 	&\multicolumn{4}{c|}{CIs for $\sigma$} \\
 \cline{4-11}
 &&&			\multicolumn{2}{c|}{$T_1$} 	&\multicolumn{2}{c|}{$T_3$}
&			\multicolumn{2}{c|}{$T_2$} 	&\multicolumn{2}{c|}{$T_4$}\\
\cline{4-11}	
$k$&$\theta$	&$n$	&AL	&CP		&AL		&CP	&AL	&CP		&AL		&CP\\

\hline
3&	0.75&	2	&20.0257 &1.0000 &20.0257 &1.0000 &31.6359 &1.0000&31.6363&1.0000\\
&&		3	&3.6177 &1.0000 &3.6177 &1.0000 &6.9250 &1.0000&6.9250&1.0000\\
&&		4	&1.7214 &1.0000 &1.7214 &1.0000 &3.9357 &1.0000&3.9357&1.0000\\
&&		5	&1.1022 &1.0000 &1.1022 &1.0000 &3.0354 &1.0000&3.0354&1.0000\\
&&		6	&0.7655 &1.0000 &0.7655 &1.0000 &2.4669 &1.0000&2.4669&1.0000\\
\hline
			\end{tabular}
	}
			\end{center}
			\end{table}

\begin{table}[H]
{\bf Table 12:} The EBs and EMSPEs of the point predictors of $R_{n+1(k)}$ and ALs and CPs of PIs for  $R_{n+1(k)}$ based on $T_1^*$ and $T_2^*$.\vspace{-0.1cm}
\begin{center}
\resizebox{\textwidth}{!}{
\begin{tabular}{|c|c|c|c|c|c|c|c|c|c|c|}
\hline
 &&&			\multicolumn{4}{c|}{Point predictors of $R_{n+1(k)}$} 	&\multicolumn{4}{c|}{PIs for $R_{n+1(k)}$} \\
 \cline{4-11}
 &&&			\multicolumn{2}{c|}{BLUP} 	&\multicolumn{2}{c|}{BLIP}
&			\multicolumn{2}{c|}{$T^*_1$} 	&\multicolumn{2}{c|}{$T^*_2$}\\
\cline{4-11}	
$k$&$\theta$	&$n$	&EB	&EMSPE		&EB		&EMSPE	&AL	&CP		&AL		&CP\\

\hline
1&	0.75&	2	&$-0.0016$ &0.0215 &0.0478 &0.0169 &4.2411 &0.9515&4.2411&0.9515\\
&&		3	&0.0021 &0.0065 &0.0190 &0.0061 &0.6387 &0.9505&0.6387&0.9505\\
&&		4	&$-0.0019$ &0.0020 &0.0059 &0.0019 &0.3027 &0.9580&0.3027&0.9580\\
&&		5	&0.0011 &0.0012 &0.0051 &0.0012 & 0.1670&0.9400&0.1670&0.9400\\
&&		6	&0.0014 &0.0005 &0.0038 &0.0005 &0.1134 &0.9455&0.1134&0.9455\\
\hline
			\end{tabular}
	}
			\end{center}
			\end{table}

\begin{table}[H]
{\bf Table 12:} continued...\vspace{-0.5cm}
\begin{center}
\resizebox{\textwidth}{!}{
\begin{tabular}{|c|c|c|c|c|c|c|c|c|c|c|}
\hline
 &&&			\multicolumn{4}{c|}{Point predictors of $R_{n+1(k)}$} 	&\multicolumn{4}{c|}{PIs for $R_{n+1(k)}$} \\
 \cline{4-11}
 &&&			\multicolumn{2}{c|}{BLUP} 	&\multicolumn{2}{c|}{BLIP}
&			\multicolumn{2}{c|}{$T^*_1$} 	&\multicolumn{2}{c|}{$T^*_2$}\\
\cline{4-11}	
$k$&$\theta$	&$n$	&EB	&EMSPE		&EB		&EMSPE	&AL	&CP		&AL		&CP\\

\hline
2&	0.75&	2	&0.0076 &0.0248 &0.0614 &0.0204 &4.0951 &0.961&4.0951&0.9610\\
&&		3	&0.0040 &0.0097 &0.0269 &0.0093 &0.7667 &0.9475&0.7667&0.9475\\
&&		4	&$-0.0008$ &0.0042 &0.0104 &0.0040 &0.3772 &0.9550&0.3772&0.9550\\
&&		5	&$-0.0012$ &0.0022 &0.0048 &0.0022 &0.2340 &0.9475&0.2340&0.9475\\
&&		6	&$-0.0011$ &0.0011 &0.0024 &0.0010 &0.1590 &0.9515&0.1590&0.9515\\
\hline
			\end{tabular}
	}
			\end{center}
			\end{table}

\begin{table}[H]
{\bf Table 12:} continued...\vspace{-0.1cm}
\begin{center}
\resizebox{\textwidth}{!}{
\begin{tabular}{|c|c|c|c|c|c|c|c|c|c|c|}
\hline
 &&&			\multicolumn{4}{c|}{Point predictors of $R_{n+1(k)}$} 	&\multicolumn{4}{c|}{PIs for $R_{n+1(k)}$} \\
 \cline{4-11}
 &&&			\multicolumn{2}{c|}{BLUP} 	&\multicolumn{2}{c|}{BLIP}
&			\multicolumn{2}{c|}{$T^*_1$} 	&\multicolumn{2}{c|}{$T^*_2$}\\
\cline{4-11}	
$k$&$\theta$	&$n$	&EB	&EMSPE		&EB		&EMSPE	&AL	&CP		&AL		&CP\\

\hline
3&	0.75&	2	&$-0.0002$ &0.0232 & 0.0561&0.0181 &3.7855 &0.9445&3.7855&0.9445\\
&&		3	&$-0.0005$ &0.0096 &0.0250 &0.0089 &0.8357 &0.9460&0.8357&0.9460\\
&&		4	&$-0.0001$ &0.0051 &0.0133 &0.0049 &0.4261 &0.9500&0.4261&0.950\\
&&		5	& $-0.0016$&0.0028 &0.0057 &0.0027 &0.2864 &0.9490&0.2864&0.9490\\
&&		6	&0.0001 &0.0015 &0.0046 &0.0015 &0.1799 &0.9425&0.1799&0.9425\\
\hline
			\end{tabular}
	}
			\end{center}
			\end{table}
From Tables 10, 11 and 12, we can reach the following conclusions.
\begin{itemize}
\item
As expected, the BLIEs possess smaller EMSEs than the corresponding BLUEs and the BLIPs possess smaller EMSPEs than the corresponding BLUPs. So the results of the simulation confirm the theoretical results  in (\ref{rec1}) and (\ref{rec2}) which state that both of the BLIEs of $\mu$ and $\sigma$ are more efficient than the corresponding BLUEs. Besides, the simulation results verify the results of Table 8, which state that the BLIP of the next future $k$-record works better than the corresponding BLUP in terms of MSPE.
\item
We observe that the EMSE decreases as the number of $k$-records increases. Moreover, we see that the EMSPE is also decreasing w.r.t. the number of $k$-records. The AL values of the CIs and PIs are also  decreasing w.r.t. the number of $k$-records. The AL values for $n=2$ are large in comparison with those for $n>2$, so we conclude that the CIs and PIs perform very poorly when the number of $k$-records equals 2.
\item
The results of CIs for $\mu$ based on $T_1$ and $T_3$ are too close to each other, implying that both $T_1$ and $T_3$ may lead to the same CIs as $T_1$ is a linear function of $T_2$. We may come to similar conclusions after checking the
results related to the CIs for $\sigma$ based on $T_2$ and $T_4$ and the results related to the PIs based on $T^*_1$ and $T_2^*$.
\end{itemize}

\section{An Illustrative Example}
In this section, we consider a real data set related to positive rates due to COVID-19 in Andorra from 30 December 2021 to 30 January 2022. The data are as follows:
\begin{verbatim}
0.2012, 0.2557, 0.2508, 0.2463, 0.2609, 0.2835, 0.3226, 0.2962,
0.3823, 0.4016, 0.4230, 0.5830, 0.6134, 0.5808, 0.5692, 0.5386,
0.5289, 0.5189, 0.2791, 0.1954, 0.1421, 0.4703, 0.4428, 0.4385,
0.4347, 0.4309, 0.7393, 0.8955, 0.5688, 0.6026, 0.7379, 0.9515,
\end{verbatim}
(see the website of COVID-19 data: https://ourwordindata.org/coronavirus-source-data).

\vspace{0.25cm}
We applied the Kolmogorov-Smirnov (K-S) test to checking if the UG distribution with pdf (1.2) fits the above data well. We obtained the K-S test statistic and its corresponding $p$-value to be $D=0.11062$ and $p=0.7883$, respectively, which emphasizes that the UG model possesses a quite good fit to the data.
From these data, we can notice that the ordinary lower records are as follows
\begin{verbatim}
0.2012,  0.1954, 0.1421.
\end{verbatim}
So we have only three lower records. However, the extracted lower 2-records are given by
\begin{verbatim}
0.2557,  0.2508, 0.2463, 0.2012, 0.1954.
\end{verbatim}
We have plotted the lower 2-records against the corresponding expected values in Table 1 for $\theta=1.5$ and we observed that there is a strong correlation between the lower 2-records and the corresponding means (correlation coefficient as high as 0.8686669), so we may conclude that the 2-record come from a UG distribution with $\theta=1.5$.

\vspace{0.25cm}
First, we estimate the location and scale parameters based on the first $n=4$ lower 2-records. From Tables 3 and 4, we obtain the BLUEs of $\mu$ and $\sigma$ as follows
\begin{eqnarray*}
\mu^* &=&a_1 R_{1(2)}+a_2 R_{2(2)}+a_3 R_{3(2)}+a_4 R_{4(2)}\\
&=&(-1.45813 \times 0.2557)+ (-0.34770 \times 0.2508)+ (-0.47173 \times 0.2463)+(3.27756 \times 0.2012)\\&=&0.08321097.
 \\&&\\
\sigma^* &=&b_1 R_{1(1)}+b_2 R_{2(1)}+b_3 R_{3(1)}+ b_4 R_{4(2)}\\
&=&(3.01524 \times 0.2557)+ (0.50480 \times 0.2508)+ (0.68667 \times 0.2463)+(-4.20671 \times 0.2012)\\&=&0.2203375.
\end{eqnarray*}
From Table 5, we compute the variances of the BLUEs and the covariance  of $\mu^*$ and $\sigma^*$ as follows
\[
Var(\mu^*)=0.07735\, \sigma^2, \quad Var(\sigma^*)= 0.19233 \,\sigma^2,  \quad  Cov(\mu^*, \sigma^*)=-0.11334\, \sigma^2 .
\]
The BLIEs of the location and scale parameters are computed to be $\widetilde{\mu}=0.1041557$ and $\widetilde{\sigma}=0.1847957$, respectively.
Moreover, we have
\[
Var(\widetilde{\mu})=0.05754\, \sigma^2, \quad Var(\widetilde{\sigma})= 0.135286 \,\sigma^2,  \quad  Cov(\widetilde{\mu}, \widetilde{\sigma})=-0.079724\, \sigma^2 .
\]
From (\ref{cimu}) and (\ref{cimu2}) and Tables 6 and 7, the $95\%$ CIs for $\mu$ based on $T_1$ and $T_3$ are given by
$(-0.3351658, 0.1506618)$ and $(-0.3351654, 0.1506602)$, respectively.
Besides, from (\ref{cisigma}) and (\ref{cisigma2}) and Tables 6 and 7, the $95\%$ CIs for $\mu$ based on $T_2$ and $T_4$ are given by
$(0.9661532, 0.1178024)$ and $(0.966052, 0.1178009)$, respectively.

\vspace{0.25cm}
Next, suppose that
we intend to predict $5^{th}$ 2-record based on the first four lower 2-records. From Table 1, we have $\alpha_{n+1}=\alpha_5=0.45806$ when $\theta$ equals 1.5.
From Table 2, we see that for $\theta=1.5$, we have
\begin{eqnarray*}
{\bm \omega}^T&=&\big(Cov(Z_{1(2)},Z_{5(2)}), Cov(Z_{2(2)},Z_{5(2)}), Cov(Z_{3(2)},Z_{5(2)}), Cov(Z_{4(2)},Z_{5(2)})\big)\\
&=&(0.00566, 0.00806, 0.00902, 0.00931).
\end{eqnarray*}
From Table 1, the vector of means of the standard 2-records for $\theta=1.5$ is given by
\[
{\bm \alpha}^T=(0.7992, 0.6680, 0.57653, 0.50939),
\]
and ${\bf B}_{4\times 4}$ is the variance-covariance matrix of the first four standard 2-record values that can be derived from Table 2. So the BLUP of $R_{5(2)}$ is computed to be $R^*_{5(2)}=0.188666$. Besides, we have
${\rm MSPE}(R^*_{5(2)})=0.001251258\, \sigma^2$.
The BLIP of $R_{5(2)}$ is calculated to be $\widetilde{R}_{5(2)}=0.1911442$ and its MSPE is computed to be
${\rm MSPE}(\widetilde{R}_{5(2)})=0.001099848\, \sigma^2$.
The real value of $R_{5(2)}$ is $0.1952$ and we see that though both $R^*_{5(2)}$ and $\widetilde{R}_{5(2)}$ are close to this value and their MSPEs seem to be small, the BLIP $R_{5(2)}$ is closer to the real value of $R_{5(2)}$ as expected.

\vspace{0.25cm}
Now, from (\ref{cimu3}), (\ref{cimu4}) and Table 9, the $95\%$ PIs for $R_{5(2)}$ based on $T^*_1$ and $T^*_2$ are given by
$(0.1235971, 0.2009576)$
and
$(0.1236043, 0.2009598)$
respectively.
As we see the real value of $R_{5(2)}$ belongs to the calculated PIs.

\section{Concluding Remarks}
In this paper, we have considered the UG distribution as the base model and have focused on the linear estimation of its location and scale parameters. We have computed the means, variances and covariances of the lower $k$-record values extracted from the standardized UG model for selected values of the shape parameter $\theta$, $n=1(1)6$ and $k=1, 2$ and 3. We then obtained the best linear unbiased and best linear invariant estimators for the location and scale parameters of the UG distribution. We have also discussed linear prediction of a future lower $k$-record value. Construction of confidence intervals for the location and scale parameters as well as prediction intervals for a future $k$-record are investigated as well. The theoretical results and the results from Table 8 emphasize that the linear invariant estimation and linear invariant prediction do work superior to the linear unbiased estimation and linear unbiased prediction, respectively, in the sense of mean squared (prediction) error criterion. This conclusion has been confirmed by the numerical results of the simulation study. Thus, we can recommend a researcher to select linear invariant estimation and prediction in the case of the UG model if s/he wants to consider the mean squared (prediction) error criterion, and of course to select linear unbiased estimation and prediction if the unbiasedness criterion is to be considered. One other point concluded from the simulation study is that the CIs for $\mu$ based on $T_1$ and $T_3$ lead to too similar numerical results, and the same point can be observed for the CIs for $\sigma$ based on $T_2$ and $T_4$ and the PIs based on $T^*_1$ and $T_2^*$. We may explain these points  by  the linear relations between the BLIEs and BLUEs and the linear relation between the BLIP  and the BLUP. Besides, the number of $k$-records affects the performance of point and interval estimators and predictors.  For example, the CIs and PIs perform very poorly when the number of $k$-records equals 2 and  we should construct CIs for the parameters and PIs for the future $k$-records based on more than 2 $k$-records.

Throughout the paper, we assume that the shape parameter $\theta$ is known and the location and scale parameters are unknown. However, in real situations, the shape parameter may not be known, so we may replace $\theta$ with its estimate, for example, its maximum likelihood or moment estimate, and then use the results of this paper. In such a case, for example, the BLUEs of the location and scale parameters are then logically called the approximate BLUEs. Another approach that may be applied to such situations is the one that is used in our real data example, namely, plotting the observed $k$-records against the expected means for several values of $\theta$ and checking if there is a strong correlation between the $k$-records and the corresponding means. Summing up, we may conclude that the results of this paper can be useful in the estimation and prediction of $k$-records in real phenomena.

All the computations of the paper were done using  the statistical software R (R Core Team, 2023), and  the packages {\sf expint} (see Goulet  {\it et al.}, 2022), {\sf essentials} (see Simmons, 2022)  and  {\sf AdequacyModel} (see Marinho \textit{et al.}, 2016)  therein.

\end{document}